\documentclass[12pt,a4paper]{article}
\usepackage[centertags]{amsmath}
\usepackage{amsfonts}
\usepackage{amssymb}
\usepackage{amsthm}
\usepackage{graphics}
\usepackage{newlfont}
\theoremstyle{plain}
\newtheorem{thm}{Theorem}[section]

\newtheorem{lem}[thm]{Lemma}
\newtheorem{prop}[thm]{Proposition}
\theoremstyle{definition}
\newtheorem{defn}[thm]{Definition}
\theoremstyle{remark}
\newtheorem{rem}[thm]{Remark}
\numberwithin{equation}{section}

\newcommand{\bq }{\begin{equation}}
\newcommand{\eq }{\end{equation}}
\newcommand{\bbb }{\begin{eqnarray}}
\newcommand{\eee }{\end{eqnarray}}
\newcommand{\bb }{\begin{eqnarray*}}
\newcommand{\ee }{\end{eqnarray*}}
\newcommand{\ed }{\end{document}}
\begin{document}
\title {Induced and Coinduced Representations of Hopf Group Coalgebra}
\author{ A.S.Hegazi$^1$, F.Ismail$^2$ and M.M. Elsofy$^3$\\
\small{$^1$ Department of Mathematics, Faculty of Science,
Mansoura University, Egypt.}\\
\small{$^2$ Department of Mathematics, Faculty of Science,
Cairo University, Egypt. \:\:\:\;\;\;\;\mbox{}    }\\
\small{$^3$ Department of Mathematics, Faculty of Science, Fayoum
University, Egypt. \;\;\;\mbox{}    }}
\date{}
\maketitle
\begin{abstract}
 In this work we study the induction theory for Hopf group coalgebra. To
reach this goal we define a substructure B of a Hopf group coalgebra $H$,
 called subHopf group coalgebra.  Also, we introduced \ the definition of
Hopf group suboalgebra and group coisotropic quantum subgroup of $H$.
\end{abstract}
\section{\large{Introduction}}

The induced representation of quantum group (quasitriangular Hopf algebra [M] and [Mon])
is introduced by Gonzalez-Ruiz, L. A. Ibort [G-I] and Ciccoli [C]. The
structure of Ciccoli has many difficultites in the structure of quantum
subgroups. For this, recently, Hegazi, Agauany, F. Ismail and I. Saleh in
[H], have succeeded to give a new algebric structure for quantum subgroups
and subquantum groups, that is a subspase $B$ of a bialgebra $H$ is called a
sub-bi-algebra if the restriction to $B$ of the structure of $H$ turns $B$
into bialgebra. But a bi-sub-algebra of $H$ is a different thing. It is pair
$(B,\pi )$ consisting of another bialgebra and surjective homomorphism $\pi
:H\rightarrow B.$ A relation between these two notions are given if the
bialgebra allows a decomposition of the form $H=B\oplus I$ with $I$ a
bi-ideal in $H$. The projection $\pi $ into B yields a bi-sub-algebra $%
(B,\pi )$ in the above sense. Conversely if $B\subseteq H$ is both a
sub-bi-algebra and a bi-sub-algebra of $H$ such that $\pi ^{2}=\pi ,$ then $%
H=B+Ker$ $\pi $ and we have a decomposition of the form $H=B\oplus I$. These
constractions made us able to introduce the induced representation of Hopf
alagebra in each case. That is a $B$-comodule for a sub-Hopf algebra $B$ of
a Hopf algebra $H$ gives rise to an induced $H$-Hopf module and that a
representation of a quantum subgroup $(B,\pi )$ of a quantum group $H$
induced a (so-called first type) Hopf representation of H. This procedure
realizes a quantum group induced representation. These structre was given
for first time by E. Tuft.\\

Recently, Turaev and   Virelizier gives use a new definition for a generalization
of Hopf algebra, for Hopf algebra structure see [S] and [Mon], i.e. Hopf group coalgebra. These generalization gives us a
new quantum group structure and the algebraic structure of these structure
is of great importance:
\begin{enumerate}
    \item  Hegazi and Abd Hafez introduced "multiplier Hopf group coalgebra".
    \item Van Deale, Hegazi and Abd Hafez introduced the new algebraic structure
"group-cograded Hopf *-algebra", J. Algebra.
    \item Hegazi, studied Differential calculus on Hopf group coalgebra (M. Sc.).
\end{enumerate}

In this paper, unless otherwise, every thing takes place over a field $K$,
and $K$-space means vector space over $K$. A map $f$ \ from a space $V$ into
a space $W$ always means linear map over K. The tensor product $V\otimes W$
is understood to be $V\otimes _{K}W$ , $I:V\rightarrow V$ always denotes the
identity map, and the transposition map $\tau :V\otimes W\rightarrow
W\otimes V$ is defined by $\tau (v\otimes w)=w\otimes v$ \ for $v\in V,w\in
W.$ Let $f:C\rightarrow D$ be a map. Then $f$ $^{\ast }:D^{\ast
}\rightarrow C^{\ast }$ is a map, where $f^{\ast }(\phi )(c)=\phi (f(c))$
for all $\phi \in D^{\ast },c\in C.$\\

Recently, quasitriangular Hopf $\pi -$coalgebra are introduced by Truaev
[T]. He gives rise to crossed $\pi $-catogeries. Virelizier [V-V1] studied the
algebraic properties of the Hopf $\pi -$coalgebras, also he has show that
the existence of integrals and trace for such catogery and has generlized
the main properties of the quasitriangular Hopf algebra to the setting of
Hopf $\pi -$coalgebra. Now, let us give some basic definitions about Hopf $%
\pi -$coalgebra. For group $\pi $, a $\pi $-coalgebra (over $K$) is a family
$C=\{C_{\alpha }\}_{\alpha \in \pi }$ of $K$-spaces endowed with a family $%
K- $maps (the comultiplication) $\Delta =\{\Delta _{\alpha ,\beta
}:H_{\alpha \beta }\rightarrow H_{\alpha }\otimes H_{\beta }\}_{\alpha
,\beta \epsilon \pi },$ and $K$-map (the counit) $\varepsilon
:H_{1}\rightarrow k$ such that the following diagrams are commute:

\begin{equation*}
\begin{tabular}{l}
$%
\begin{tabular}{lll}
$\quad C_{\alpha \beta \gamma }$ &
\begin{tabular}{l}
$\Delta _{\alpha ,\beta \gamma }$ \\
$\longrightarrow $%
\end{tabular}
& $C_{\alpha }\otimes C_{\beta \gamma }$ \\
$\Delta _{\alpha \beta ,\gamma }\downarrow $ &  & $\downarrow I_{\alpha
}\otimes \Delta _{\beta ,\gamma }$ \\
$\quad C_{\alpha \beta }\otimes C_{\gamma }$ &
\begin{tabular}{l}
$\longrightarrow $ \\
$\Delta _{\alpha ,\beta }\otimes I_{\gamma }$%
\end{tabular}
& $C_{\alpha }\otimes C_{\beta }\otimes C_{\gamma }$%
\end{tabular}%
\ $ \\
 \mbox{}\\
$\quad (\Delta _{\alpha ,\beta }\otimes I_{\gamma }$ $)\Delta _{\alpha \beta
,\gamma }=(I_{\alpha }\otimes \Delta _{\beta ,\gamma })\Delta _{\alpha
,\beta \gamma }$%
\end{tabular}%
\end{equation*}

\begin{equation*}
\begin{tabular}{l}
\begin{tabular}{lll}
& $C_{1}\otimes C_{\alpha }$ &  \\
$\epsilon \otimes I_{\alpha }\swarrow $ &  &  \\
 \mbox{}\\
$\quad K\otimes C_{\alpha }$ & $\quad \uparrow \Delta _{1,\alpha }$ &  \\
 \mbox{}\\
$\qquad \quad \sim _{\alpha }\nwarrow $ &  &  \\
& $\quad C_{\alpha }$ &
\end{tabular}
$\quad $\quad $\quad $\quad
\begin{tabular}{ll}
$C_{\alpha }\otimes C_{1}$ &  \\
& $\searrow I_{\alpha }\otimes \epsilon $ \\
 \mbox{}\\
$\Delta _{\alpha ,1}\uparrow $ & $C_{\alpha }\otimes K$ \\
 \mbox{}\\
& $\nearrow \sim _{\alpha }$ \\
$\qquad C_{\alpha }$ &  \\
&
\end{tabular}
$\quad $\quad $\quad $\quad $\quad $\quad  \\
$\quad $\quad $\ \ \ \ \ \ \ \ \ \quad \ \ \ \ \ \ \ \ (\epsilon \otimes
I_{\alpha })\Delta _{1,\alpha }=$ $\sim _{\alpha }=(I_{\alpha }\otimes
\epsilon )\Delta _{\alpha ,1}$%
\end{tabular}%
\ .
\end{equation*}

A Hopf $\pi -$coalgebra is a $\pi $-coalgebra $(H=\{H_{\alpha }\}_{\alpha
\in \pi },\Delta ,\epsilon )$ with a family $S=\{S_{\alpha }:H_{\alpha
}\rightarrow H_{\alpha ^{-1}}\}_{\alpha \in \pi }$ of $K$-maps such that

\begin{enumerate}
\item $H_{\alpha }$ is an algebra with multiplication $\mu _{\alpha }$ and
unit $\eta _{\alpha }(1_{K})$ for all $\alpha \in \pi ;$

\item $\Delta _{\alpha ,\beta }$ $,$ $\epsilon $ are algebra maps for all $%
\alpha ,\beta \in \pi ,$

\item The antipode $S$ satisfy
\begin{equation*}
\begin{tabular}{l}
\begin{tabular}{lllll}
$\quad \qquad H_{1}$ &  & $\quad \overset{1_{\alpha }\epsilon }{\rightarrow }
$ &  & $H_{\alpha }$ \\
$\Delta _{\alpha ^{-1},\alpha }\downarrow $ &  &  &  & $\uparrow \mu
_{\alpha }$ \\
$H_{\alpha ^{-1}}\otimes H_{\alpha }$ &  & $\overset{S_{\alpha ^{-1}}\otimes
I_{\alpha }}{\longrightarrow }$ &  & $H_{\alpha }\otimes H_{\alpha }$%
\end{tabular}
\ \ \ \ \ \
\begin{tabular}{lll}
$\quad \qquad H_{1}$ & $\quad \overset{1_{\alpha }\epsilon }{\rightarrow }$
& $H_{\alpha }$ \\
$\Delta _{\alpha ,\alpha ^{-1}}\downarrow $ &  & $\uparrow \mu _{\alpha }$
\\
$H_{\alpha }\otimes H_{\alpha ^{-1}}$ & $\overset{I_{\alpha }\otimes
S_{\alpha ^{-1}}}{\longrightarrow }$ & $H_{\alpha }\otimes H_{\alpha }$%
\end{tabular}
\\
 \mbox{}\\
$\qquad \qquad \mu _{\alpha }(S_{\alpha ^{-1}}\otimes I_{\alpha })\Delta
_{\alpha ^{-1},\alpha }=1_{\alpha }\epsilon =\mu _{\alpha }(I_{\alpha
}\otimes S_{\alpha ^{-1}})\Delta _{\alpha ,\alpha ^{-1}}$%
\end{tabular}%
\
\end{equation*}
\end{enumerate}

A Hopf $\pi $-coalgebra $H$ is of finite type when every $H_{\alpha }$ is
finite dimensional. Note that it does not means that ${\oplus _{\alpha \in
\pi }}H_{\alpha }$ is finite dimensional (unless $H_{\alpha }=0$ for all but
a finite number of $\alpha \in \pi $). $H$ is totally finite when the direct
sum ${\oplus _{\alpha \in \pi }}H_{\alpha }$ is finite dimensional. The
notion of Hopf $\pi -$coalgebra is not self dual, i.e., given a Hopf $\pi -$%
coalgebra $H=\{H_{\alpha }\}_{\alpha \in \pi }$, the family $H^{\ast
}=\{H_{\alpha }^{\ast }\}_{\alpha \in \pi }$ does not have a natural
structure of Hopf $\pi -$coalgebra. Note, $(H_{1},\Delta _{1,1},\epsilon )$
is a (classical) Hopf algebra.\newline

A Hopf $\pi -$coalgebra $H$ is said to have a left (right) cosection if
there exist a family of algebra maps $\eta =\{\eta _{\alpha
}:H_{1}\rightarrow H_{\alpha }\}_{\alpha \in \pi }$ such that $\eta _{\alpha
}\,$is left (right) $H_{1}$-comodule map for all $\alpha \in \pi $ i.e., the
following digram is commute
\begin{equation*}
\begin{tabular}{llll}
& $H_{1}$ & $\overset{\eta _{\alpha }}{\longrightarrow }$ & $H_{\alpha }$ \\
$\Delta _{1,1}$ & $\downarrow $ &  & $\downarrow \Delta _{1,\alpha }$ \\
&  &  &  \\
& $H_{1}\otimes H_{1}$ &
\begin{tabular}{l}
$\longrightarrow $ \\
$I_{1}\otimes \eta _{\alpha }$%
\end{tabular}
& $H_{1}\otimes H_{\alpha }$%
\end{tabular}%
\end{equation*}%

Let $C$ be a $\pi $-coalgebra. A right $\pi $-comodule over $C$ is a family $%
M=\{M_{\alpha }\}_{\alpha \in \pi }$ of $K$-spaces endowed with a family $%
\theta =\{\theta _{\alpha ,\beta }:M_{\alpha \beta }\rightarrow M_{\alpha
}\otimes C_{\beta }\}_{\alpha ,\beta \epsilon \pi }$ of $K$-maps such that
the following diagrams are commute:
\begin{equation*}
\begin{tabular}{lll}
$%
\begin{tabular}{lll}
$\qquad M_{\alpha \beta \gamma }$ & $\quad \overset{\theta _{\alpha ,\beta
\gamma }}{\longrightarrow }$ & $\quad M_{\alpha }\otimes C_{\beta \gamma }$
\\
$\theta _{\alpha \beta ,\gamma }\downarrow $ &  & $\qquad \downarrow
I_{\alpha }\otimes \Delta _{\beta ,\gamma }$ \\
$\quad M_{\alpha \beta }\otimes C_{\gamma }$ & \quad $\overset{\theta
_{\alpha ,\beta }\otimes I_{\gamma }}{\longrightarrow }$ & $C_{\alpha
}\otimes C_{\beta }\otimes C\gamma $%
\end{tabular}%
\ $ &  &
\begin{tabular}{l}
\end{tabular}
\begin{tabular}{ll}
$\quad M_{\alpha }$ &  \\
& $\searrow \theta _{\alpha ,1}$ \\
$\quad \sim _{M_{\alpha }}\downarrow $ & $\quad M_{\alpha }\otimes C_{1}$ \\
& $\swarrow I_{\alpha }\otimes \epsilon $ \\
$M_{\alpha }\otimes K$ &
\end{tabular}
\\
\qquad \qquad \qquad \qquad \qquad \qquad \qquad \qquad &  &  \\
$\quad (\theta _{\alpha ,\beta }\otimes I_{\gamma }$ $)\theta _{\alpha \beta
,\gamma }=(I_{\alpha }\otimes \Delta _{\beta ,\gamma })\theta _{\alpha
,\beta \gamma }$ &  & $\qquad (I_{\alpha }\otimes \epsilon )\theta _{\alpha
,1}=\sim _{M_{\alpha }}$%
\end{tabular}%
\end{equation*}%
\newline

A right Hopf $\pi $-comodule over $H$ is a right $\pi $-comodule $%
M=\{M_{\alpha }\}_{\alpha \in \pi }$ by coaction $\theta =\{\theta _{\alpha
,\beta }\}_{\alpha ,\beta \epsilon \pi }$ such that

\begin{enumerate}
\item $M_{\alpha }$ is $H_{\alpha }-$module by action $\rho _{\alpha }$ for
all $\alpha \in \pi $

\item The following diagram, for all $\alpha ,\beta ,\gamma \in \pi $, is
commute:%
\begin{equation*}
\begin{tabular}{l}
\begin{tabular}{llllll}
$\ \ \ \quad M_{\alpha \beta }\otimes H_{\alpha \beta }$ & $\overset{\rho
_{\alpha \beta }}{\rightarrow }$ & $\quad M_{\alpha \beta }$ & $\overset{%
\theta _{\alpha ,\beta }}{\rightarrow }$ & $M_{\alpha }\otimes H_{\beta }$ &
\\
$\theta _{\alpha ,\beta }\otimes \Delta _{\alpha ,\beta }\downarrow $ &  &
&  & $\ \ \ \ \ \uparrow \rho _{\alpha }\otimes \mu _{\beta }$ &  \\
$\ \ \ \ (M_{\alpha }\otimes H_{\beta })\otimes (H_{\alpha }\otimes H_{\beta
})$ &  & $\overset{I_{\alpha }\otimes \tau \otimes I_{\beta }}{\rightarrow }$
&  & $(M_{\alpha }\otimes H_{\alpha })\otimes (H_{\beta }\otimes H_{\beta })$
&
\end{tabular}
\\
\begin{tabular}{l}
\\
$\qquad i.e.,$ $\theta _{\alpha ,\beta }$ $\rho _{\alpha \beta }=(\rho
_{\alpha }\otimes \mu _{\beta })(I_{\alpha }\otimes \tau \otimes I_{\beta
})(\theta _{\alpha ,\beta }\otimes \Delta _{\alpha ,\beta })$%
\end{tabular}%
\end{tabular}%
\end{equation*}%
\newline
\end{enumerate}

Let C be a $\pi -$coa$\lg $ebra. A $\pi -$coideal of $C$ is a family $%
V=\{V_{\alpha }\}_{\alpha \in \pi }$ such \nolinebreak that

\begin{enumerate}
\item $V_{\alpha }$ is subspace of $C_{\alpha }$ for all $\alpha \in \pi ,$

\item $\Delta _{\alpha ,\beta }(V_{\alpha \beta })\subset V_{\alpha }\otimes
C_{\beta }+C_{\alpha }\otimes V_{\beta }$ for all $\alpha ,\beta \in \pi .$

\item $\epsilon _{1}(V_{1})=0.$
\end{enumerate}

A Hopf $\pi $-coideal of $H$ is a family of $K$-spaces $V=\{V_{\alpha
}\}_{\alpha \in \pi }$ such that $V$ is $\pi $-coideal of $H$, $V_{\alpha }$
is ideal of $H_{\alpha }$ and $S_{\alpha }(V_{\alpha })\subseteq V_{\alpha
^{-1}}$ for all $\alpha \in \pi .$ \newline

A family $A=\{A_{\alpha }\}_{\alpha \in \pi }$ is called subHopf $\pi -$coa$%
\lg $ebra of $H$ $\ $if $(A,\Delta ,\epsilon )$ is $\pi -$coalgebra of $H$, $%
A_{\alpha }$ is a subalgebra of $H_{\alpha }{}$ and $S_{\alpha }(A_{\alpha
})\subseteq A_{\alpha ^{-1}}$for all $\alpha \in \pi .$

A family $A=\{A_{\alpha }\}_{\alpha \in \pi }$ of a subalgebra of Hopf $\pi $%
-coalgebra of $H$ is called isolated if there exist a family $I=\{I_{\alpha
}\}_{\alpha \in \pi }$ of Hopf $\pi $- coideal of $H$ such that$H_{\alpha
}=A_{\alpha }\oplus I_{\alpha }$for all $\alpha \in \pi .$\newline

We will call Hopf $\pi -$subcoalgebra of $H$ any pair $(C,\sigma )$ such
that $C=\{C_{\alpha }\}_{\alpha \in \pi }$ is a Hopf $\pi -$coalgebra and $%
\sigma =\{\sigma _{\alpha }:H_{\alpha }\rightarrow C_{\alpha }\}_{\alpha \in
\pi }$ family of algebra epimorphisms which satisfies for $\alpha ,\beta \in
\pi $

\begin{enumerate}
\item $\Delta _{\alpha ,\beta }^{C}\sigma _{\alpha \beta }=(\sigma _{\alpha
}\otimes \sigma _{\beta })\Delta _{\alpha ,\beta }^{H},$ i.e., the following
digram is commute
\begin{equation*}
\begin{tabular}{llll}
& $H_{\alpha \beta }$ & $\overset{\sigma _{\alpha ,\beta }}{\longrightarrow }
$ & $C_{\alpha \beta }$ \\
$\Delta _{\alpha ,\beta }^{H}$ & $\downarrow $ &  & $\downarrow \Delta
_{\alpha ,\beta }^{C}$ \\
&  &  &  \\
& $H_{\alpha }\otimes H_{\beta }$ &
\begin{tabular}{l}
$\longrightarrow $ \\
$\sigma _{\alpha }\otimes \sigma _{\beta }$%
\end{tabular}
& $C_{\alpha }\otimes C_{\beta }$%
\end{tabular}%
\end{equation*}

\item $\varepsilon ^{C} \sigma _{1}=\varepsilon ^{H},$ i.e., the following
digram is commute
\begin{equation*}
\begin{tabular}{lllll}
& $H_{1}$ & $\overset{\sigma _{1}}{\longrightarrow }$ & $C_{1}$ &  \\
$\epsilon ^{H}$ & $\searrow $ &  & $\swarrow $ & $\epsilon ^{C}$ \\
&  & $K$ &  &
\end{tabular}%
\end{equation*}

\item $S_{\alpha }^{C}\sigma _{\alpha }=\sigma _{\alpha ^{-1}}S_{\alpha
}^{H} $ ($\sigma _{\beta \alpha \beta ^{-1}}\varphi _{\beta }^{H}{}=\varphi
_{\beta }^{C}\sigma _{\alpha }$) i.e., the following digram is commute
\begin{equation*}
\begin{tabular}{lllll}
& $H_{\alpha }$ & $\overset{\sigma _{\alpha }}{\longrightarrow }$ & $%
C_{\alpha }$ &  \\
$S_{\alpha }^{H}$ & $\downarrow $ &  & $\downarrow $ & $S_{\alpha }^{C}$ \\
& $H_{\alpha ^{-1}}$ &
\begin{tabular}{l}
$\longrightarrow $ \\
$\sigma _{\alpha ^{-1}}$%
\end{tabular}
& $C_{\alpha ^{-1}}$ &
\end{tabular}%
\end{equation*}
\end{enumerate}

A pair $(C,\sigma )$ is called left $\pi -$coisotropic quantum subgroup of $%
H $ if

\begin{enumerate}
\item $C$ is $\pi -$coa$\lg $ebra$,$

\item $C_{\alpha }$ is left $H_{\alpha }$-module by $\omega _{\alpha }$ for
all $\alpha \in \pi ,$

\item $\sigma =\{\sigma _{\alpha }:H_{\alpha }\rightarrow C_{\alpha }$ $%
\}_{\alpha \in \pi }$ family of surjective linear maps such that

\begin{enumerate}
\item $\sigma _{\alpha }$ is left $H_{\alpha }$-module map for all $\alpha
\in \pi ,$

\item $\Delta _{\alpha ,\beta }^{C}\sigma _{\alpha \beta }=(\sigma _{\alpha
}\otimes \sigma _{\beta })\Delta _{\alpha ,\beta }^{H},$

\item $\varepsilon ^{C}\circ \sigma _{1}=\varepsilon ^{H}.$\newline
\end{enumerate}
\end{enumerate}

\begin{prop} Every isolated subHopf $\pi -$coalgebra
is $\pi -$coisotropic quantum subgroup.
\end{prop}

\begin{proof} Clear from definition.
\end{proof}

A left $\pi -$coisotropic quantum subgroup $(C,\sigma )$ of $H$ is said to
have a left section if there exist a family of linear, convolution
invertible, maps $g=\{g_{\alpha }:C_{\alpha }\rightarrow H_{\alpha
}\}_{\alpha \in \pi }$ such that

\begin{enumerate}
\item $g_{\alpha }(\sigma _{\alpha }(1))=1$

\item For $\alpha \in \pi ,u\in H_{1},c\in C_{\alpha }$ and $v\in \sigma
_{\alpha }^{-1}(c)$%
\begin{eqnarray*}
&&(\sigma _{1}\otimes I_{\alpha }^{H})(\mu _{1}\otimes I_{\alpha
})(I_{1}\otimes \tau )(\Delta _{1,\alpha }^{H}g_{\alpha }\otimes
I_{1})(c\otimes u) \\
&=&(I_{1}\otimes g_{\alpha })(\sigma _{1}\otimes \sigma _{\alpha })(\mu
_{1}\otimes I_{\alpha })(I_{1}\otimes \tau )(\Delta _{1,\alpha }^{H}\otimes
I_{1})(v\otimes u)
\end{eqnarray*}
\end{enumerate}

\section{\large{Induced representations of Hopf
group coalgebra}}

In this section, we study the induced representation for Hopf group
coalgebra. To reach this goal we use the definitions of subHopf group
coalgebra, Hopf group subcoalgebra and group coisotropic quantum subgroup.

\begin{thm}Let $H=\{H_{\alpha }\}_{\alpha \in \pi }$ be
a Hopf $\pi -$coalgebra and $(C,\sigma )$ be a Hopf $\pi -$subcoalgebra of $%
H.$ Then $(C,\sigma )$ is a left $\pi -$coisotropic quantum subgroup of $H.$
\end{thm}
\begin{proof}
We define $\omega =\{\omega _{\alpha }=\mu _{\alpha
}^{C}(\sigma _{\alpha }\otimes I_{\alpha }):H_{\alpha }\otimes C_{\alpha
}\rightarrow C_{\alpha }\}_{\alpha \in \pi }$ we`ll prove that $C_{\alpha }$
is a left $H_{\alpha }-$module by $\rho _{\alpha }.$ i.e., the followig
digrams are commute
\begin{equation*}
\begin{tabular}{llll}
& $H_{\alpha }\otimes H_{\alpha }\otimes C_{\alpha }$ & $%
\begin{tabular}{l}
$I_{\alpha }\otimes \omega _{\alpha }$ \\
$\ \ \rightarrow $%
\end{tabular}%
$ & $H_{\alpha }\otimes C_{\alpha }$ \\
&  &  &  \\
$\mu _{\alpha }\otimes I_{\alpha }$ & $\downarrow $ &  & $\downarrow \omega
_{\alpha }$ \\
& $H_{\alpha }\otimes C_{\alpha }$ &
\begin{tabular}{l}
$\longrightarrow $ \\
$\rho _{\alpha }$%
\end{tabular}
& $C_{\alpha }$%
\end{tabular}%
\mbox{}\hspace{.5cm}\mbox{}\\
\begin{tabular}{lll}
$K\otimes C_{\alpha }$ & $%
\begin{tabular}{l}
$\eta _{\alpha }\otimes I_{\alpha }$ \\
$\ \ \rightarrow $%
\end{tabular}%
$ & $H_{\alpha }\otimes C_{\alpha }$ \\
$\sim \searrow $ &  & $\swarrow \omega _{\alpha }$ \\
& $C_{\alpha }$ &
\end{tabular}%
\end{equation*}

\begin{eqnarray*}
\omega _{\alpha }(I_{\alpha }\otimes \omega _{\alpha })(h\otimes k\otimes
b)
&=&\omega _{\alpha }(h\otimes \sigma _{\alpha }(k)b) \\
&=&\sigma _{\alpha }(h)(\sigma _{\alpha }(k)b) \\
&=&\sigma _{\alpha }(hk)b \\
&=&\omega _{\alpha }(hk\otimes b) \\
&=&\omega _{\alpha }(\mu _{\alpha }\otimes I_{\alpha })(h\otimes k\otimes b)
\end{eqnarray*}

and%
\begin{eqnarray*}
\omega _{\alpha }(\eta _{\alpha }\otimes I_{\alpha })(k\otimes b)
&=&\omega _{\alpha }(\eta _{\alpha }(k)\otimes b) \\
&=&\sigma _{\alpha }(\eta _{\alpha }(k))b \\
&=&kb \\
&=&\sim (k\otimes b)
\end{eqnarray*}

Now, we`ll prove that $\sigma _{\alpha }$ is left $H_{\alpha }$-module map
for all $\alpha \in \pi $ i.e., the following digram is commute
\begin{equation*}
\begin{tabular}{llll}
& $H_{\alpha }$ & $%
\begin{tabular}{l}
$\sigma _{\alpha }$ \\
$\rightarrow $%
\end{tabular}%
$ & $C_{\alpha }$ \\
&  &  &  \\
$\mu _{\alpha }^{H}$ & $\uparrow $ &  & $\uparrow \omega _{\alpha }$ \\
& $H_{\alpha }\otimes H_{\alpha }$ &
\begin{tabular}{l}
$\longrightarrow $ \\
$I_{\alpha }\otimes \sigma _{\alpha }$%
\end{tabular}
& $H_{\alpha }\otimes C_{\alpha }$%
\end{tabular}%
\end{equation*}

\begin{equation*}
\omega _{\alpha }(I_{\alpha }\otimes \sigma _{\alpha })=\mu _{\alpha
}^{C}(\sigma _{\alpha }\otimes I_{\alpha })(I_{\alpha }\otimes \sigma
_{\alpha })=\mu _{\alpha }^{C}(\sigma _{\alpha }\otimes \sigma _{\alpha
})=\sigma _{\alpha }\mu _{\alpha }^{H}.
\end{equation*}%
\end{proof}

\begin{rem}
If $(C,\sigma )$ is a Hopf $\pi -$subcoalgebra of $H$,
then the map
\begin{equation*}
L_{\alpha ,\beta }=(\sigma _{\alpha }\otimes I_{\beta }^{H})\Delta _{\alpha
,\beta }^{H}:H_{\alpha \beta }\rightarrow C_{\alpha }\otimes H_{\beta
}~~\forall \alpha ,\beta \in \pi
\end{equation*}%
is an algebra map.
\end{rem}
\begin{prop}
Let $H=\{H_{\alpha }\}_{\alpha \in \pi }$ be a Hopf
$\pi -$coalgebra and $(C,\sigma )$ be a Hopf $\pi -$subcoalgebra of $H.$ Let
$B=\{B_{\alpha }\}_{\alpha \in \pi },$ where
$
B_{\alpha }=\{h\in H_{\alpha }:L_{1,\alpha }(h)=1\otimes h\},
$
then $B_{\alpha }$ is subalgebra of $H_{\alpha }$ and $B$ is right $\pi -$%
coideal of $H$ (this B is called $\pi -$quantum right embeddable homogeneous
space of $H$ ).
\end{prop}

\begin{proof}\mbox

\begin{enumerate}
\item Let $h,k\in B_{\alpha },$ then
\begin{equation*}
L_{1,\alpha }(hk)=L_{1,\alpha }(h)\cdot L_{1,\alpha }(k)=(1\otimes h)\cdot
(1\otimes k)=1\otimes hk
\end{equation*}
and hence $hk\in B_{\alpha }.$

\item Let $h\in B_{\alpha \beta },$ we`ll prove that $\Delta _{\alpha ,\beta
}^{H}(h)\in B_{\alpha }\otimes H_{\beta }.$
\begin{eqnarray*}
(L_{1,\alpha }\otimes I_{\beta }^{H})\Delta _{\alpha ,\beta }^{H}(h)
&=&((\sigma _{1}\otimes I_{\alpha }^{H})\Delta _{1,\alpha }^{H}\otimes
I_{\beta }^{H})\Delta _{\alpha ,\beta }^{H}(h) \\
&=&(\sigma _{1}\otimes I_{\alpha }^{H}\otimes I_{\beta }^{H})(\Delta
_{1,\alpha }^{H}\otimes I_{\beta }^{H})\Delta _{\alpha ,\beta }^{H}(h) \\
&=&(\sigma _{1}\otimes I_{\alpha }^{H}\otimes I_{\beta
}^{H})(I_{1}^{H}\otimes \Delta _{\alpha ,\beta }^{H})\Delta _{1,\alpha \beta
}^{H}(h) \\
&=&(I_{1}^{C}\otimes \Delta _{\alpha ,\beta }^{H})(\sigma _{1}\otimes
I_{\alpha \beta }^{H})\Delta _{1,\alpha \beta }^{H}(h) \\
&=&(I_{1}^{C}\otimes \Delta _{\alpha ,\beta }^{H})(1\otimes h) \\
&=&1\otimes \Delta _{1,\alpha \beta }^{H}(h)
\end{eqnarray*}
\end{enumerate}
\end{proof}

\begin{lem}
Let $H=\{H_{\alpha }\}_{\alpha \in \pi }$ be a Hopf $\pi
- $coalgebra and let $(C,\sigma )$ be a left $\pi -$coisotropic quantum
subgroup of $H.$ For $\alpha ,\beta ,\gamma \in \pi $ and $a,b\in H_{\alpha
\beta }$ we have

\begin{itemize}
\item $L_{\alpha ,\beta }(ab)=\Delta _{\alpha ,\beta }^{H}(a)\Theta
L_{\alpha ,\beta }(b)$

\item $(I_{\alpha }\otimes \Delta _{\beta ,\gamma }^{H})L_{\alpha ,\beta
\gamma }=(L_{\alpha ,\beta }\otimes I_{\gamma })\Delta _{\alpha \beta
,\gamma }^{H}$,
\end{itemize}
 where $(m\otimes n)\Theta (u\otimes v)=\omega _{\alpha }(m\otimes u)\otimes
nv,m,u\in H_{\alpha },n,v\in H_{\beta }.$

\end{lem}

\begin{proof}

\begin{eqnarray*}
L_{\alpha ,\beta }(ab)&=&(\sigma _{\alpha }\otimes I_{\beta }^{H})\Delta
_{\alpha ,\beta }^{H}\mu _{\alpha \beta }(a\otimes b) \\
&=&(\sigma _{\alpha }\otimes I_{\beta }^{H})(\mu _{\alpha }\otimes \mu
_{\beta })(I\otimes \tau \otimes I)(\Delta _{\alpha ,\beta }^{H}\otimes
\Delta _{\alpha ,\beta }^{H})(a\otimes b)\\
&=&(\sigma _{\alpha }\mu _{\alpha }\otimes \mu _{\beta })(I\otimes \tau
\otimes I)(\Delta _{\alpha ,\beta }^{H}\otimes \Delta _{\alpha ,\beta
}^{H})(a\otimes b) \\
&=&(\sigma _{\alpha }\mu _{\alpha }\otimes \mu _{\beta })(I\otimes \tau
\otimes I)(\Delta _{\alpha ,\beta }^{H}\otimes \Delta _{\alpha ,\beta
}^{H})(a\otimes b) \\
&=&(\omega _{\alpha }(I_{\alpha }\otimes \sigma _{\alpha })\otimes \mu
_{\beta })(a_{1}^{\alpha }\otimes a_{2}^{\beta }\otimes b_{1}^{\alpha
}\otimes b_{2}^{\beta }) \\
&=&\omega _{\alpha }(a_{1}^{\alpha }\otimes \sigma _{\alpha }(b_{1}^{\alpha
}))\otimes a_{2}^{\beta }b_{2}^{\beta } \\
&=&\Delta _{\alpha ,\beta }^{H}(a)\Theta L_{\alpha ,\beta }(b)
\end{eqnarray*}%
Also,%
\begin{eqnarray*}
(I_{\alpha }\otimes \Delta _{\beta ,\gamma }^{H})L_{\alpha ,\beta \gamma }
&=&(I_{\alpha }\otimes \Delta _{\beta ,\gamma }^{H})(\sigma _{\alpha
}\otimes I_{\beta \gamma }^{H})\Delta _{\alpha ,\beta \gamma }^{H} \\
&=&(\sigma _{\alpha }\otimes I_{\beta }^{H}\otimes I_{\gamma
}^{H})(I_{\alpha }\otimes \Delta _{\beta ,\gamma }^{H})\Delta _{\alpha
,\beta \gamma }^{H} \\
&=&(\sigma _{\alpha }\otimes I_{\beta }^{H}\otimes I_{\gamma }^{H})(\Delta
_{\alpha ,\beta }^{H}\otimes I_{\gamma })\Delta _{\alpha \beta ,\gamma }^{H}
\\
&=&((\sigma _{\alpha }\otimes I_{\beta }^{H})\Delta _{\alpha ,\beta
}^{H}\otimes I_{\gamma }^{H})\Delta _{\alpha \beta ,\gamma }^{H} \\
&=&(L_{\alpha ,\beta }\otimes I_{\gamma })\Delta _{\alpha \beta ,\gamma }^{H}
\end{eqnarray*}
\end{proof}

\begin{prop}
Let $H=\{H_{\alpha }\}_{\alpha \in \pi }$ be a Hopf
$\pi -$coalgebra and let $(C,\sigma )$ be a left $\pi -$coisotropic quantum
subgroup of $H.$ Let $G=\{G_{\alpha }\}_{\alpha \in \pi },$ where $G_{\alpha
}=\{h\in H_{\alpha }:L_{1,\alpha }(h)=\sigma _{1}(1)\otimes h\}$, then $%
G_{\alpha }$ is subalgebra of $H_{\alpha }$ and $G$ is right $\pi -$ coideal
of $H$.
\end{prop}

\begin{proof}
For $h,g\in G_{\alpha },$ we have

\begin{enumerate}
\item $L_{1,\alpha }(h)=\sigma _{1}(1)\otimes h$ and $L_{1,\alpha
}(g)=\sigma _{1}(1)\otimes g$

\item $\omega _{1}(I_{1}\otimes \sigma _{1})=\sigma _{1}\mu _{1}$
\end{enumerate}
We`ll prove that $hg\in G_{\alpha },i.e.,L_{1,\alpha }(hg)=\sigma
_{1}(1)\otimes hg$
\begin{eqnarray*}
L_{1,\alpha }(hg)
&=& \Delta_{1,\alpha}(h)\Theta L_{1,\alpha}(g)\;\;\;\;\;\;\;\;\;\;\;\;\;\;\;\;\;\ \textrm{by Lemma 2.4} \\
&=& \Delta_{1,\alpha}(h)\Theta (\sigma_{1}(1)\otimes g)\\
&=&(\omega _{1}\otimes \mu _{\alpha })(h_{1}^{1}\otimes \sigma
_{1}(1)\otimes h_{2}^{\alpha }\otimes g) \\
&=&\omega _{1}(I_{1}\otimes \sigma _{1})(h_{1}^{1}\otimes 1)\otimes
h_{2}^{\alpha }g \\
&=&\sigma _{1}\mu _{1}(h_{1}^{1}\otimes 1)\otimes h_{2}^{\alpha }g \\
&=&\sigma _{1}(h_{1}^{1})\otimes h_{2}^{\alpha }g \\
&=&(\sigma _{1}\otimes \mu _{\alpha })(\Delta _{1,\alpha }^{H}\otimes
I_{\alpha })(h\otimes g) \\
&=&(I_{1}\otimes \mu _{\alpha })((\sigma _{1}\otimes I_{1})\Delta _{1,\alpha
}^{H}(h)\otimes g) \\
&=&(I_{1}\otimes \mu _{\alpha })(\sigma _{1}(1)\otimes h\otimes g) \\
&=&\sigma _{1}(1)\otimes hg
\end{eqnarray*}
Now, we`ll prove that $G$ is right $\pi $-coideal of $H$, i.e., for $h\in
G_{\alpha \beta },\Delta _{\alpha ,\beta }^{H}(h)\in G_{\alpha }\otimes
H_{\beta }$
\begin{eqnarray*}
(L_{1,\alpha }\otimes I_{\beta }^{H})\Delta _{\alpha ,\beta }^{H}(h)
&=&((\sigma _{1}\otimes I_{\alpha }^{H})\Delta _{1,\alpha }^{H}\otimes
I_{\beta }^{H})\Delta _{\alpha ,\beta }^{H}(h) \\
&=&(\sigma _{1}\otimes I_{\alpha }^{H}\otimes I_{\beta }^{H})(\Delta
_{1,\alpha }^{H}\otimes I_{\beta }^{H})\Delta _{\alpha ,\beta }^{H}(h) \\
&=&(\sigma _{1}\otimes I_{\alpha }^{H}\otimes I_{\beta
}^{H})(I_{1}^{H}\otimes \Delta _{\alpha ,\beta }^{H})\Delta _{1,\alpha \beta
}^{H}(h) \\
&=&(I_{1}^{C}\otimes \Delta _{\alpha ,\beta }^{H})(\sigma _{1}\otimes
I_{\alpha \beta }^{H})\Delta _{1,\alpha \beta }^{H}(h) \\
&=&(I_{1}^{C}\otimes \Delta _{\alpha ,\beta }^{H})(\sigma _{1}(1)\otimes h)
\\
&=&\sigma _{1}(1)\otimes \Delta _{1,\alpha \beta }^{H}(h)
\end{eqnarray*}%
\end{proof}

\begin{thm}
Suppose $V=\{V_{\alpha }\}_{\alpha \in \pi },$ be a
right $\pi -$comodule over the left $\pi -$coisotropic quantum subgroup $C$
of Hopf $\pi -$coalgebra $H$ by $\rho =\{\rho _{\alpha ,\beta }:V_{\alpha
\beta }\rightarrow V_{\alpha }\otimes C_{\beta }\}_{\alpha ,\beta \in \pi }.$
Then $Ind(\rho )=\{Ind(\rho )_{\alpha }\}_{\alpha \in \pi }$ where
$
Ind(\rho )_{\alpha }=\{x\in V_{1}\otimes H_{\alpha }:(I_{1}\otimes
L_{1,\alpha })x=(\rho _{1,1}\otimes I_{\alpha })x\}
$
is right $\pi -$comodule over $H$ by $(I\otimes \Delta )=\{(I\otimes \Delta
)_{\alpha ,\beta }=I_{1}\otimes \Delta _{\alpha ,\beta }\}_{\alpha ,\beta
\in \pi }$
\end{thm}

\begin{proof}
 We will prove that, for $\alpha ,\beta  \in \pi ,$
\begin{equation*}
(I\otimes \Delta )_{\alpha ,\beta }(Ind(\rho )_{\alpha \beta })\subseteq
Ind(\rho )_{\alpha }\otimes H_{\beta }.
\end{equation*}

Since, for $v\otimes h\in Ind(\rho )_{\alpha \beta },$ we have
\begin{eqnarray*}
(I_{1}\otimes L_{1,\alpha }\otimes I_{\beta })(I_{1}\otimes \Delta
_{\alpha ,\beta })(v\otimes h)
&=&(I_{1}\otimes (L_{1,\alpha }\otimes I_{\beta })\Delta _{\alpha ,\beta
})(v\otimes h) \\
&=&(I_{1}\otimes (I_{1}\otimes \Delta _{\alpha ,\beta })L_{1,\alpha \beta
})(v\otimes h)\;\;\;\;\;\;\;\;\; \textrm{by Lemma 2.4  }\\
&=&(I_{1}\otimes I_{1}\otimes \Delta _{\alpha ,\beta })(I_{1}\otimes
L_{1,\alpha \beta })(v\otimes h) \\
&=&(I_{1}\otimes I_{1}\otimes \Delta _{\alpha ,\beta })(\rho _{1,1}\otimes
I_{\alpha \beta })(v\otimes h) \\
&=&(\rho _{1,1}\otimes I_{\alpha }\otimes I_{\beta })(I_{1}\otimes \Delta
_{\alpha ,\beta })(v\otimes h)
\end{eqnarray*}

Now, we`ll prove that the following digrams are commute
\begin{equation*}
\begin{tabular}{llll}
& $Ind(\rho )_{\alpha \beta \gamma }$ & $\overset{(I\otimes \Delta )_{\alpha
\beta ,\gamma }}{\longrightarrow }$ & $Ind(\rho )_{\alpha \beta }\otimes
H_{\gamma }$ \\
&  &  &  \\
$(I\otimes \Delta )_{\alpha ,\beta \gamma }$ & $\downarrow $ &  & $%
\downarrow (I\otimes \Delta )_{\alpha ,\beta }\otimes I_{\gamma }$ \\
&  &  &  \\
& $Ind(\rho )_{\alpha }\otimes H_{\beta \gamma }$ &
\begin{tabular}{l}
$\longrightarrow $ \\
$I_{Ind(\rho )_{\alpha }}\otimes \Delta _{\beta ,\gamma }$%
\end{tabular}
& $Ind(\rho )_{\alpha }\otimes H_{\beta }\otimes H_{\gamma }$%
\end{tabular}%
\end{equation*}%
\mbox{and}\\
\begin{equation*}
\begin{tabular}{llll}
$Ind(\rho )_{\alpha }$ &  & $%
\begin{tabular}{l}
$(I\otimes \Delta )_{\alpha ,1}$ \\
$\rightarrow $%
\end{tabular}%
$ & $Ind(\rho )_{\alpha }\otimes H_{1}$ \\
&  &  &  \\
& $\sim _{Ind(\rho )_{\alpha }}\searrow $ &  & $\swarrow I_{Ind(\rho
)_{\alpha }}\otimes \epsilon ^{H}$ \\
&  & $Ind(\rho )_{\alpha }\otimes K$ &
\end{tabular}%
\end{equation*}

\begin{eqnarray*}
((I\otimes \Delta )_{\alpha ,\beta }\otimes I_{\gamma })(I\otimes \Delta
)_{\alpha \beta ,\gamma } &=&(I_{1}\otimes \Delta _{\alpha ,\beta }\otimes I_{\gamma })(I_{1}\otimes \Delta _{\alpha \beta ,\gamma }) \\
&=&(I_{1}\otimes (\Delta _{\alpha ,\beta }\otimes I_{\gamma })\Delta
_{\alpha \beta ,\gamma }) \\
&=&(I_{1}\otimes (I_{\alpha }\otimes \Delta _{\beta ,\gamma })\Delta
_{\alpha ,\beta \gamma }) \\
&=&(I_{1}\otimes I_\alpha\otimes\Delta_{\beta,\gamma})(I_1\otimes \Delta_{\alpha,\beta\gamma})\\
&=& (I_{Ind(\rho )_{\alpha }}\otimes \Delta_{\beta,\gamma})(I\otimes\Delta)_{\alpha,\beta\gamma}
\end{eqnarray*}

\begin{eqnarray*}
(I_{Ind(\rho )_{\alpha }}\otimes \epsilon ^{H})(I\otimes \Delta )_{\alpha
,1}&=&(I_{1}\otimes I_{\alpha }\otimes \epsilon ^{H})(I_{1}\otimes \Delta
_{\alpha ,1}) \\
&=&(I_{1}\otimes (I_{\alpha }\otimes \epsilon )\Delta _{\alpha
,1})\\
&=&(I_{1}\otimes \sim _{H_{\alpha }})~=\sim _{Ind(\rho )_{\alpha }}
\end{eqnarray*}
\end{proof}

\begin{rem}
Given a right corepresentation $\rho $ of the $\pi -$%
coisotropic quantum subgroup of $(C,\sigma )$ the corresponding
corepresentation $(I\otimes \Delta)$ on $Ind(\rho )$ of $H$ is called induced representation from $\rho$ on $H$.
\end{rem}

\section{\large{ Geometric realization for induced representation}}

Throughout this section $H$ is Hopf $\pi -$ coalgebra, $(C,\sigma )$ is Hopf
$\pi -$subcoalgebra of $H$ and $V=\{V_{\alpha }\}_{\alpha \in \pi }$ be a
right $\pi -$comodule over $C$. The purpose of this section is to explicit
that the induced representation $Ind(\rho )$ from Hopf group subcoalgebra $H$
is isomorphic to the tensor product of $\pi -$quantum embeddable homogeneous
space $B$ (in Proposition 2.3) with the given comodule $V$ as module and in
case $(C,\sigma )$ is left $\pi -$coisotropic quantum subgroup $H$ is
isomorphic to $C\otimes G$ as vector space where $G=\{G_{\alpha }\}_{\alpha
\in \pi },$ where
\begin{equation*}
G_{\alpha }=\{h\in H_{\alpha }:L_{1,\alpha }(h)=(\sigma _{1}\otimes
I_{\alpha }^{H})\Delta _{1,\alpha }^{H}(h)=\sigma _{1}(1)\otimes h\}.
\end{equation*}

\begin{lem}
$Ind(\rho
)_{\alpha }$ is a right $B_{\alpha }$-module for all $\alpha \in \pi .$
\end{lem}

\begin{proof}
Let $v_{1}\otimes h_\alpha \in Ind(\rho)_\alpha, b_{\alpha} \in B_{\alpha},$ we define the right action as follows
\begin{equation*}
\lambda _{\alpha }(v_{1}\otimes h_{\alpha }\otimes b_{\alpha })=v_{1}\otimes
h_{\alpha }b_{\alpha }
\end{equation*}
We need only to prove that
\begin{equation*}
\lambda _{\alpha }(v_{1}\otimes h_{\alpha }\otimes b_{\alpha })=v_{1}\otimes
b_{\alpha }h_{\alpha }\in Ind(\rho )_{\alpha }.
\end{equation*}
We have%
\begin{eqnarray*}
(I_{1}\otimes L_{1,\alpha })(v_{1}\otimes h_{\alpha })
&=&(\rho _{1,1}\otimes I\smallskip _{\alpha })(v_{1}\otimes h_{\alpha })
\end{eqnarray*}%
imply that%
\begin{eqnarray*}
((I_{1}\otimes L_{1,\alpha })(v_{1}\otimes h_{\alpha }))(1\otimes 1\otimes
b_{\alpha })
&=&((\rho _{1,1}\otimes I_{\alpha })(v_{1}\otimes h_{\alpha }))(1\otimes
1\otimes b_{\alpha })
\end{eqnarray*}%
imply that%
\begin{eqnarray}\label{2}
v_{1}\otimes \sigma _{1}((h_{\alpha })_{1})\otimes (h_{\alpha
})_{2}b_{\alpha }
&=&\rho _{1,1}(v_{1})\otimes h_{\alpha }b_{\alpha }.
\end{eqnarray}
Also, we have
\begin{eqnarray} \label{1}
L_{1,\alpha }(h_{\alpha }b_{\alpha })
&=&\Delta _{1,\alpha }^{H}(h_{\alpha })\Theta L_{1,\alpha }(b_{\alpha }) \;\;\;\;\;\textrm{by Lemma 2.5}\nonumber\\
&=&((h_{\alpha })_{1}\otimes (h_{\alpha })_{2})\Theta (1_{C_{1}}\otimes
b_{\alpha }) \nonumber\\
&=&\omega _{\alpha }((h_{\alpha })_{1}\otimes \sigma _{1}(1_{H_{1}}))\otimes
(h_{\alpha })_{2}b_{\alpha } \nonumber\\
&=&\sigma _{1}(h_{\alpha _{1}})\otimes (h_{\alpha })_{2}b_{\alpha }.
\end{eqnarray}%
Now,
\begin{eqnarray*}
(I_{1}\otimes L_{1,\alpha })\lambda _{\alpha }(v_{1}\otimes h_{\alpha
}\otimes b_{\alpha })
&=&(I_{1}\otimes L_{1,\alpha })(v_{1}\otimes h_{\alpha }b_{\alpha }) \\
&=&v_{1}\otimes L_{1,\alpha }(h_{\alpha }b_{\alpha }) \\
&=&v_{1}\otimes \sigma _{1}((h_{\alpha })_{1})\otimes (h_{\alpha
})_{2}b_{\alpha } \;\;\;\;\;\;\;\textrm{by Equation \ref{1}}\\
&=&\rho _{1,1}(v_{1})\otimes h_{\alpha }b_{\alpha }\;\;\;\;\;\;\;\;\;\;\;\;\;\;\;\;\;\;\;\;\;\;\textrm{by Equation \ref{2}} \\
&=&(\rho _{1,1}\otimes I_{\alpha })(v_{1}\otimes h_{\alpha
}b_{\alpha })\\
&=& (\rho _{1,1}\otimes I_{\alpha }) \lambda _{\alpha }(v_{1}\otimes h_{\alpha
}\otimes b_{\alpha }).
\end{eqnarray*}
The left action is similar.
\end{proof}

\begin{defn}
Let $(C,\sigma )$ be a Hopf $\pi -$subcoalgebra of $%
H.$ If $g=\{g_{\alpha }:C_{\alpha }\rightarrow H_{\alpha }\}_{\alpha \in \pi
},$ is a family of linear maps, its convolution inverse, if it exists, is a family of
linear maps $g^{-1}=\{g_{\alpha }^{-1}:C_{\alpha ^{-1}}\rightarrow H_{\alpha
}\}_{\alpha \in \pi }$ such that
\begin{equation*}
\mu _{\alpha }(g_{\alpha }\otimes g_{\alpha }^{-1})\Delta _{\alpha ,\alpha
^{-1}}^{C}(c)=\epsilon ^{C}(c)1_{H_{\alpha }}=\mu _{\alpha }(g_{\alpha
}^{-1}\otimes g_{\alpha })\Delta _{\alpha ^{-1},\alpha }^{C}(c)
\end{equation*}%
\end{defn}

\begin{defn}
A Hopf $\pi -$subcoalgebra $(C,\sigma )$ of $H$ is
said to have a left section if there exists a family of linear, convolution
invertible, maps $g=\{g_{\alpha }:C_{\alpha }\rightarrow H_{\alpha
}\}_{\alpha \in \pi }$ such that for all $\alpha\in \pi$,

\begin{enumerate}
\item $g_{\alpha }(1)=1$

\item $L_{1,\alpha }g_{\alpha }=(I_{1}\otimes g_{\alpha })\Delta _{1,\alpha
}^{C}$
\end{enumerate}
\end{defn}

\begin{lem}
For any Hopf $\pi -$coalgebra $H$ we have
\begin{enumerate}
\item $(\Delta _{1,\alpha }^{H}\otimes \Delta _{\alpha ^{-1},1}^{H})\Delta
_{\alpha ,\alpha ^{-1}}^{H}=(I_{1}\otimes \Delta _{\alpha ,\alpha
^{-1}}^{H}\otimes I_{1})(\Delta _{1,1}^{H}\otimes I_{1})\Delta _{1,1}^{H}$

\item $(\Delta _{\alpha ^{-1},1}^{H}\otimes \Delta _{1,\alpha }^{H})\Delta
_{\alpha ^{-1},\alpha }^{H}=(I_{\alpha ^{-1}}\otimes \Delta
_{1,1}^{H}\otimes I_{\alpha })$($\Delta _{\alpha ^{-1},1}^{H}\otimes
I_{\alpha })\Delta _{\alpha ^{-1},\alpha }^{H}\newline
$
\end{enumerate}
\end{lem}

\begin{proof}
\begin{eqnarray*}
(\Delta _{1,\alpha }^{H}\otimes \Delta _{\alpha ^{-1},1}^{H})\Delta
_{\alpha ,\alpha ^{-1}}^{H}
&=&(\Delta _{1,\alpha }^{H}\otimes I_{\alpha ^{-1}}\otimes I_{1})(I_{\alpha
}\otimes \Delta _{\alpha ^{-1},1}^{H})\Delta _{\alpha ,\alpha ^{-1}}^{H}%
\newline
\\
&=&(\Delta _{1,\alpha }^{H}\otimes I_{\alpha ^{-1}}\otimes I_{1})(\Delta
_{\alpha ,\alpha ^{-1}}^{H}\otimes I_{1})\Delta _{1,1}^{H} \\
&=&((\Delta _{1,\alpha }^{H}\otimes I_{\alpha ^{-1}})\Delta _{\alpha ,\alpha
^{-1}}^{H}\otimes I_{1})\Delta _{1,1}^{H} \\
&=&((I_{1}\otimes \Delta _{\alpha ,\alpha ^{-1}}^{H})\Delta
_{1,1}^{H}\otimes I_{1})\Delta _{1,1}^{H} \\
&=&(I_{1}\otimes \Delta _{\alpha ,\alpha ^{-1}}^{H}\otimes I_{1})(\Delta
_{1,1}^{H}\otimes I_{1})\Delta _{1,1}^{H}
\end{eqnarray*}%
\begin{equation*}
i.e.,(\Delta _{1,\alpha }^{H}\otimes \Delta _{\alpha ^{-1},1}^{H})\Delta
_{\alpha ,\alpha ^{-1}}^{H}(h)=h_{11}^{1}\otimes h_{121}^{\alpha }\otimes
h_{122}^{\alpha ^{-1}}\otimes h_{2}^{1}
\end{equation*}%
The second statement is similar.
\end{proof}

Now, from Lemma \textbf{\ }3.5\textbf{\ }up to Theorem 3.12 below $(C,\sigma
)$ have a section $g=\{g_{\alpha }\}_{\alpha \in \pi }$ and antipode $S^C=\{S^{C}_{\alpha}: C_\alpha \to C_{\alpha^{-1}}\}_{\alpha\in\pi}$

\begin{lem}
For $\alpha \in \pi ,L_{1,\alpha }g_{\alpha
}^{-1}=(S_{1}^{C}\otimes g_{\alpha }^{-1})\tau \Delta _{\alpha ^{-1},1}^{C}$

\end{lem}

\begin{proof}
We`ll prove that for $\alpha \in \pi $, $L_{1,\alpha
}g_{\alpha }^{-1}$ and $(S_{1}^{C}\otimes g_{\alpha }^{-1})\tau \Delta
_{\alpha ^{-1},1}^{C}$ are inverse to the same element $L_{1,\alpha
}g_{\alpha }$ in the convolution algebra $Conv(C,C_{1}\otimes H_{\alpha }).$%
\begin{eqnarray*}
(L_{1,\alpha }g_{\alpha }\ast L_{1,\alpha }g_{\alpha }^{-1})(c)
&=&\{(\sigma _{1}\otimes I_{\alpha }^{H})\Delta _{1,\alpha }^{H}g_{\alpha
}\ast (\sigma _{1}\otimes I_{\alpha }^{H})\Delta _{1,\alpha }^{H}g_{\alpha
}^{-1}\}(c) \\
&=&\mu _{C_{1}\otimes H_{\alpha }}\{(\sigma _{1}\otimes I_{\alpha
}^{H})\Delta _{1,\alpha }^{H}g_{\alpha }\otimes (\sigma _{1}\otimes
I_{\alpha }^{H})\Delta _{1,\alpha }^{H}g_{\alpha }^{-1}\}\Delta _{\alpha
,\alpha ^{-1}}^{C}(c) \\
&=&(\mu _{C_{1}}\otimes \mu _{H_{\alpha }})(I\otimes \tau \otimes
I)\{(\sigma _{1}\otimes I_{\alpha }^{H})\Delta _{1,\alpha }^{H}g_{\alpha
}\otimes (\sigma _{1}\otimes I_{\alpha }^{H})\Delta _{1,\alpha
}^{H}g_{\alpha }^{-1}\}\Delta _{\alpha ,\alpha ^{-1}}^{C}(c) \\
&=&(\mu _{C_{1}}\otimes \mu _{H_{\alpha }})(\sigma _{1}\otimes \sigma
_{1}\otimes I_{\alpha }^{H}\otimes I_{\alpha }^{H}) \\
&&(I\otimes \tau \otimes I)(\Delta _{1,\alpha }^{H}\otimes \Delta _{1,\alpha
}^{H})(g_{\alpha }\otimes g_{\alpha }^{-1})\Delta _{\alpha ,\alpha
^{-1}}^{C}(c) \\
&=&(\mu _{C_{1}}(\sigma _{1}\otimes \sigma _{1})\otimes \mu _{H_{\alpha
}})(I\otimes \tau \otimes I)(\Delta _{1,\alpha }^{H}\otimes \Delta
_{1,\alpha }^{H})(g_{\alpha }\otimes g_{\alpha }^{-1})\Delta _{\alpha
,\alpha ^{-1}}^{C}(c) \\
&=& (\sigma_1 \mu_{H_1}\otimes  \mu _{H_{\alpha
}})(I\otimes \tau \otimes I)(\Delta _{1,\alpha }^{H}\otimes \Delta
_{1,\alpha }^{H})(g_{\alpha }\otimes g_{\alpha }^{-1})\Delta _{\alpha
,\alpha ^{-1}}^{C}(c) \\
&=&(\sigma _{1}\otimes I_{\alpha }^{H})(\mu _{H_{1}}\otimes \mu _{H_{\alpha
}})(I\otimes \tau \otimes I)(\Delta _{1,\alpha }^{H}\otimes \Delta
_{1,\alpha }^{H})(g_{\alpha }\otimes g_{\alpha }^{-1})\Delta _{\alpha
,\alpha ^{-1}}^{C}(c) \\
&=&(\sigma _{1}\otimes I_{\alpha }^{H})\Delta _{1,\alpha }^{H}\mu
_{H_{\alpha }}(g_{\alpha }\otimes g_{\alpha }^{-1})\Delta _{\alpha ,\alpha
^{-1}}^{C}(c) \\
&=&(\sigma _{1}\otimes I_{\alpha }^{H})\Delta _{1,\alpha }^{H}(\epsilon
^{C}(c)1_{H_{\alpha }}) \\
&=&\epsilon ^{C}(c)(\sigma _{1}\otimes I_{\alpha }^{H})\Delta _{1,\alpha
}^{H}(1_{H_{\alpha }}) \\
&=&\epsilon ^{C}(c)(1_{C_{1}}\otimes 1_{H_{\alpha }}).
\end{eqnarray*}%
and%
\begin{eqnarray*}
\{L_{1,\alpha }g_{\alpha }\ast (S_{1}^{C}\otimes g_{\alpha }^{-1})\tau
\Delta _{\alpha ^{-1},1}^{C}\}(c)
&=&\{(\sigma _{1}\otimes I_{\alpha }^{H})\Delta _{1,\alpha }^{H}g_{\alpha
}\ast (S_{1}^{C}\otimes g_{\alpha }^{-1})\tau \Delta _{\alpha
^{-1},1}^{C}\}(c) \\
&=&\{(I_{1}\otimes g_{\alpha })\Delta _{1,\alpha }^{C}\ast (S_{1}^{C}\otimes
g_{\alpha }^{-1})\tau \Delta _{\alpha ^{-1},1}^{C}\}(c) \\
&=&\mu_{C_1\otimes H_{\alpha}}\{(I_{1}\otimes g_{\alpha })\Delta _{1,\alpha }^{C}\otimes(S_{1}^{C}\otimes g_{\alpha }^{-1})\tau \Delta _{\alpha ^{-1},1}^{C}\}\Delta
_{\alpha ,\alpha ^{-1}}^{C}(c) \\
&=&(\mu _{C_{1}}\otimes \mu _{H_{\alpha }})(I\otimes \tau \otimes
I)\{(I_{1}\otimes g_{\alpha })\Delta _{1,\alpha }^{C}\otimes\\
&&(S_{1}^{C}\otimes g_{\alpha }^{-1})\tau \Delta _{\alpha ^{-1},1}^{C}\}\Delta
_{\alpha ,\alpha ^{-1}}^{C}(c) \\
&=&(\mu _{C_{1}}\otimes \mu _{H_{\alpha }})(I\otimes \tau \otimes
I)(I_{1}\otimes g_{\alpha }\otimes S_{1}^{C}\otimes g_{\alpha }^{-1}) \\
&&(I_{1}\otimes I_{\alpha }\otimes \tau )(\Delta _{1,\alpha }^{C}\otimes
\Delta _{\alpha ^{-1},1}^{C})\Delta _{\alpha ,\alpha ^{-1}}^{C}(c) \\
&=&(\mu _{C_{1}}\otimes \mu _{H_{\alpha }})(I\otimes \tau \otimes
I)(I_{1}\otimes g_{\alpha }\otimes S_{1}^{C}\otimes g_{\alpha
}^{-1})(I_{1}\otimes I_{\alpha }\otimes \tau ) \\
&&(c_{11}^{1}\otimes c_{121}^{\alpha }\otimes c_{122}^{\alpha ^{-1}}\otimes
c_{2}^{1})  \;\;\;\;\;\;\;\;\textrm{by Lemma 3.4 (1)}\\
&=&c_{11}^{1}S_{1}^{C}(c_{2}^{1})\otimes g_{\alpha }(c_{121}^{\alpha
})g_{\alpha }^{-1}(c_{122}^{\alpha ^{-1}}) \\
&=&c_{11}^{1}S_{1}^{C}(c_{2}^{1})\otimes \epsilon ^{C}(c_{12}^{\alpha
})1_{H_{\alpha }} \\
&=&\epsilon ^{C}(c_{12}^{\alpha })c_{11}^{1}S_{1}^{C}(c_{2}^{1})\otimes
1_{H_{\alpha }} \\
&=&c_{1}^{1}S_{1}^{C}(c_{2}^{1})\otimes 1_{H_{\alpha }} \\
&=&\epsilon ^{C}(c)(1_{C_{1}}\otimes 1_{H_{\alpha }}).
\end{eqnarray*}
\end{proof}

\begin{lem}
For $h\in H_{\alpha },$ we have
\begin{equation*}
h=g_{\alpha }\sigma _{\alpha }(h_{1}^{\alpha })g_{\alpha }^{-1}\sigma
_{\alpha ^{-1}}(h_{21}^{\alpha ^{-1}})h_{22}^{\alpha }.
\end{equation*}%
\end{lem}

\begin{proof}
\begin{eqnarray*}
h &=&\epsilon ^{H}(h_{1}^{1})h_{2}^{\alpha } \\
&=&\varepsilon ^{C}\sigma _{1}(h_{1}^{1})h_{2}^{\alpha } \\
&=&(\mu _{\alpha }(g_{\alpha }\otimes g_{\alpha }^{-1})\Delta _{\alpha
,\alpha ^{-1}}^{C}(\sigma _{1}(h_{1}^{1})))h_{2}^{\alpha } \\
&=&(\mu _{\alpha }(g_{\alpha }\otimes g_{\alpha }^{-1})(\sigma _{\alpha
}\otimes \sigma _{\alpha ^{-1}})\Delta _{\alpha ,\alpha
^{-1}}^{H}(h_{1}^{1}))h_{2}^{\alpha } \\
&=&\mu _{\alpha }(\mu _{\alpha }(g_{\alpha }\otimes g_{\alpha }^{-1})(\sigma
_{\alpha }\otimes \sigma _{\alpha ^{-1}})\otimes I_{\alpha })(\Delta
_{\alpha ,\alpha ^{-1}}^{H}\otimes I_{\alpha })\Delta _{1,\alpha }^{H}(h) \\
&=&\mu _{\alpha }(\mu _{\alpha }(g_{\alpha }\otimes g_{\alpha }^{-1})(\sigma
_{\alpha }\otimes \sigma _{\alpha ^{-1}})\otimes I_{\alpha })(I_{\alpha
}\otimes \Delta _{\alpha ^{-1},\alpha }^{H})\Delta _{\alpha ,1}^{H}(h) \\
&=&g_{\alpha }\sigma _{\alpha }(h_{1}^{\alpha })g_{\alpha }^{-1}\sigma
_{\alpha ^{-1}}(h_{21}^{\alpha ^{-1}})h_{22}^{\alpha }.
\end{eqnarray*}
\end{proof}

\begin{lem}
For $h\in H_{1},\alpha \in \pi ,$ we have
\begin{equation*}
\mu _{\alpha }(g_{\alpha }^{-1}\sigma _{\alpha ^{-1}}\otimes I_{\alpha
})\Delta _{\alpha ,\alpha ^{-1}}^{H}(h)\in B_{\alpha }.\newline
\end{equation*}
\end{lem}

\begin{proof}
\begin{eqnarray*}
L_{1,\alpha }\mu _{\alpha }(g_{\alpha }^{-1}\sigma _{\alpha ^{-1}}\otimes
I_{\alpha })\Delta _{\alpha ^{-1},\alpha }^{H}(h)
&=&(\sigma _{1}\otimes I_{\alpha }^{H})\Delta _{1,\alpha }^{H}\mu _{\alpha
}(g_{\alpha }^{-1}\sigma _{\alpha ^{-1}}\otimes I_{\alpha })\Delta _{\alpha
^{-1},\alpha }^{H}(h) \\
&=&\mu _{C_{1}\otimes H_{\alpha }}((\sigma _{1}\otimes I_{\alpha
}^{H})\Delta _{1,\alpha }^{H}g_{\alpha }^{-1}\sigma _{\alpha ^{-1}}\otimes
(\sigma _{1}\otimes I_{\alpha }^{H})\Delta _{1,\alpha }^{H})\Delta _{\alpha
^{-1},\alpha }^{H}(h) \\
&=&\mu _{C_{1}\otimes H_{\alpha }}((S_{1}^{C}\otimes g_{\alpha }^{-1})\tau
\Delta _{\alpha ^{-1},1}^{C}\sigma _{\alpha ^{-1}}\otimes (\sigma
_{1}\otimes I_{\alpha }^{H})\Delta _{1,\alpha }^{H})\Delta _{\alpha
^{-1},\alpha }^{H}(h) \\
&=&\mu _{C_{1}\otimes H_{\alpha }}((S_{1}^{C}\otimes g_{\alpha }^{-1})\tau
(\sigma _{\alpha ^{-1}}\otimes \sigma _{1})\Delta _{\alpha
^{-1},1}^{H}\otimes\\&& (\sigma _{1}\otimes I_{\alpha }^{H})\Delta _{1,\alpha
}^{H})\Delta _{\alpha
^{-1},\alpha }^{H}(h) \\
&=&\mu _{C_{1}\otimes H_{\alpha }}((S_{1}^{C}\otimes g_{\alpha }^{-1})\tau
\otimes I_{1}\otimes I_{\alpha }^{H}) \\
&&(\sigma _{\alpha ^{-1}}\otimes \sigma _{1}\otimes \sigma _{1}\otimes
I_{\alpha }^{H})(\Delta _{\alpha ^{-1},1}^{H}\otimes \Delta _{1,\alpha
}^{H})\Delta _{\alpha
^{-1},\alpha }^{H}(h) \\
&=&\mu _{C_{1}\otimes H_{\alpha }}((S_{1}^{C}\otimes g_{\alpha }^{-1})\tau
\otimes I_{1}\otimes I_{\alpha }^{H}) \\
&&(\sigma _{\alpha ^{-1}}(h_{11}^{\alpha ^{-1}})\otimes \sigma
_{1}(h_{121}^{1})\otimes \sigma _{1}(h_{122}^{1})\otimes h_{2}^{\alpha }) \\
&=&S_{1}^{C}(\sigma _{1}(h_{121}^{1}))\sigma _{1}(h_{122}^{1})\otimes
g_{\alpha }^{-1}\sigma _{\alpha ^{-1}}(h_{11}^{\alpha ^{-1}})h_{2}^{\alpha }
\\
&=&\epsilon ^{C}(\sigma _{1}(h_{12}^{1})1\otimes g_{\alpha }^{-1}\sigma
_{\alpha ^{-1}}(h_{11}^{\alpha ^{-1}})h_{2}^{\alpha } \\
&=&\epsilon ^{H}(h_{12}^{1})1\otimes g_{\alpha }^{-1}\sigma _{\alpha
^{-1}}(h_{11}^{\alpha ^{-1}})h_{2}^{\alpha } \\
&=&1\otimes g_{\alpha }^{-1}\sigma _{\alpha ^{-1}}(\epsilon
^{H}(h_{12}^{1})h_{11}^{\alpha ^{-1}})h_{2}^{\alpha } \\
&=&1\otimes g_{\alpha }^{-1}\sigma _{\alpha ^{-1}}(h_{1}^{\alpha
^{-1}})h_{2}^{\alpha } \\
&=&1\otimes \mu _{\alpha }(g_{\alpha }^{-1}\sigma _{\alpha ^{-1}}\otimes
I_{\alpha })\Delta _{\alpha ^{-1},\alpha }^{H}(h).
\end{eqnarray*}
\end{proof}

\begin{lem}
For $\alpha \in \pi ,$%
\begin{equation*}
(\sigma _{\alpha }\otimes \sigma _{\alpha ^{-1}}\otimes I_{\alpha
}^{H})(\Delta _{\alpha ,\alpha ^{-1}}^{H}\otimes I_{\alpha }^{H})\Delta
_{1,\alpha }^{H}g_{\alpha }=(I_{\alpha }\otimes I_{\alpha ^{-1}}\otimes
g_{\alpha })(I_{\alpha }\otimes \Delta _{\alpha ^{-1},\alpha }^{C})\Delta
_{\alpha ,1}^{C}
\end{equation*}
\end{lem}

\begin{proof}
\begin{eqnarray*}
(\sigma _{\alpha }\otimes \sigma _{\alpha ^{-1}}\otimes I_{\alpha
}^{H})(\Delta _{\alpha ,\alpha ^{-1}}^{H}\otimes I_{\alpha }^{H})\Delta
_{1,\alpha }^{H}g_{\alpha }
&=&((\sigma _{\alpha }\otimes \sigma _{\alpha ^{-1}})\Delta _{\alpha ,\alpha
^{-1}}^{H}\otimes I_{\alpha }^{H})\Delta _{1,\alpha }^{H}g_{\alpha } \\
&=&(\Delta _{\alpha ,\alpha ^{-1}}^{C}\sigma _{1}\otimes I_{\alpha
}^{H})\Delta _{1,\alpha }^{H}g_{\alpha } \\
&=&(\Delta _{\alpha ,\alpha ^{-1}}^{C}\otimes I_{\alpha }^{H})(\sigma
_{1}\otimes I_{\alpha }^{H})\Delta _{1,\alpha }^{H}g_{\alpha } \\
&=&(\Delta _{\alpha ,\alpha ^{-1}}^{C}\otimes I_{\alpha }^{H})(I_{1}\otimes
g_{\alpha })\Delta _{1,\alpha }^{C} \\
&=&(I_{\alpha }\otimes I_{\alpha ^{-1}}\otimes g_{\alpha })(\Delta _{\alpha
,\alpha ^{-1}}^{C}\otimes I_{\alpha })\Delta _{1,\alpha }^{C} \\
&=&(I_{\alpha }\otimes I_{\alpha ^{-1}}\otimes g_{\alpha })(I_{\alpha
}\otimes \Delta _{\alpha ^{-1},\alpha }^{C})\Delta _{\alpha ,1}^{C}
\end{eqnarray*}
\end{proof}

\begin{lem}
For $\alpha \in \pi ,b\in B_{\alpha }$%
\begin{equation*}
(\sigma _{\alpha }\otimes \sigma _{\alpha ^{-1}}\otimes I_{\alpha
}^{H})(\Delta _{\alpha ,\alpha ^{-1}}^{H}\otimes I_{\alpha }^{H})\Delta
_{1,\alpha }^{H}(b)=1\otimes 1\otimes b
\end{equation*}
\end{lem}

\begin{proof}
\begin{eqnarray*}
(\sigma _{\alpha }\otimes \sigma _{\alpha ^{-1}}\otimes I_{\alpha
}^{H})(\Delta _{\alpha ,\alpha ^{-1}}^{H}\otimes I_{\alpha }^{H})\Delta
_{1,\alpha }^{H}(b)
&=&((\sigma _{\alpha }\otimes \sigma _{\alpha ^{-1}})\Delta _{\alpha ,\alpha
^{-1}}^{H}\otimes I_{\alpha }^{H})\Delta _{1,\alpha }^{H}(b) \\
&=&(\Delta _{\alpha ,\alpha ^{-1}}^{C}\sigma _{1}\otimes I_{\alpha
}^{H})\Delta _{1,\alpha }^{H}(b) \\
&=&(\Delta _{\alpha ,\alpha ^{-1}}^{C}\otimes I_{\alpha }^{H})(\sigma
_{1}\otimes I_{\alpha }^{H})\Delta _{1,\alpha }^{H}(b) \\
&=&(\Delta _{\alpha ,\alpha ^{-1}}^{C}\otimes I_{\alpha }^{H})(1\otimes
b)=1\otimes 1\otimes b
\end{eqnarray*}
\end{proof}
Now, we can prove the main theorems in this section.
\begin{thm}
$H$ is isomorphic to $C\otimes B$ as vector space.
\end{thm}

\begin{proof}
We define $A=\{A_{\alpha }:C_{\alpha }\otimes B_{\alpha
}\rightarrow H_{\alpha }\}_{\alpha \in \pi }$ as follow
\begin{equation*}
A_{\alpha }(c\otimes b)=\mu _{\alpha }(g_{\alpha }\otimes I_{\alpha
})(c\otimes b)=g_{\alpha }(c)b.
\end{equation*}%
Clear $A_{\alpha }$ is linear and by Lemma 3.6 and Lemma 3.7 that $A_{\alpha
}$ is surjective for all $\alpha \in \pi .$ We define $A_{\alpha
}^{-1}=(\sigma _{\alpha }\otimes \mu _{\alpha }(g_{\alpha }^{-1}\sigma
_{\alpha ^{-1}}\otimes I_{\alpha }))(\Delta _{\alpha ,\alpha
^{-1}}^{H}\otimes I_{\alpha })\Delta _{1,\alpha }^{H}.$ We`ll prove that for
$\alpha \in \pi ,A_{\alpha }A_{\alpha }^{-1}=I_{H_{\alpha }}$ and $A_{\alpha
}^{-1}A_{\alpha }=I_{C_{\alpha }\otimes B_{\alpha }}.$ Firstly,  let $h\in H_{\alpha }$
and $c\otimes b\in C_{\alpha }\otimes B_{\alpha }$%
\begin{eqnarray*}
A_{\alpha }A_{\alpha }^{-1}(h)
&=&A_{\alpha }((\sigma _{\alpha }\otimes \mu _{\alpha }(g_{\alpha
}^{-1}\sigma _{\alpha ^{-1}}\otimes I_{\alpha }))(\Delta _{\alpha ,\alpha
^{-1}}^{H}\otimes I_{\alpha })\Delta _{1,\alpha }^{H}(h)) \\
&=&A_{\alpha }((\sigma _{\alpha }\otimes \mu _{\alpha }(g_{\alpha
}^{-1}\sigma _{\alpha ^{-1}}\otimes I_{\alpha }))(I_{\alpha }\otimes \Delta
_{\alpha ^{-1},\alpha }^{H})\Delta _{\alpha ,1}^{H}(h)) \\
&=&A_{\alpha }(\sigma _{\alpha }(h_{1}^{\alpha })\otimes g_{\alpha
}^{-1}\sigma _{\alpha ^{-1}}(h_{21}^{\alpha ^{-1}})h_{22}^{\alpha }) \\
&=&g_{\alpha }\sigma _{\alpha }(h_{1}^{\alpha })g_{\alpha }^{-1}\sigma
_{\alpha ^{-1}}(h_{21}^{\alpha ^{-1}})h_{22}^{\alpha } \\
&=&h~\ \ \ \ \ \ \ \ \ \ \ \ \ \ \ \ \ \ \ \ \ \ \ \ \ \ \ \ \ \ \ \ \ \ \ \
\ \ \ \ \ \text{by Lemma 3.6}
\end{eqnarray*}%
Secondly, since, for $\alpha \in \pi $ we have $((\sigma _{\alpha }\otimes \sigma
_{\alpha ^{-1}})\Delta _{\alpha ,\alpha ^{-1}}^{H}\otimes I_{\alpha
}^{H})\Delta _{1,\alpha }^{H}$ is an algebra map, then
\begin{eqnarray*}
A_{\alpha }^{-1}A_{\alpha }(c\otimes b)
&=&(\sigma _{\alpha }\otimes \mu _{\alpha }(g_{\alpha }^{-1}\sigma _{\alpha
^{-1}}\otimes I_{\alpha }))(\Delta _{\alpha ,\alpha ^{-1}}^{H}\otimes
I_{\alpha })\Delta _{1,\alpha }^{H}\mu _{\alpha }(g_{\alpha }\otimes
I_{\alpha })(c\otimes b) \\
&=&(I_{\alpha }\otimes \mu _{\alpha })(I_{\alpha }\otimes g_{\alpha
}^{-1}\otimes I_{\alpha })(\sigma _{\alpha }\otimes \sigma _{\alpha
^{-1}}\otimes I_{\alpha }^{H})(\Delta _{\alpha ,\alpha ^{-1}}^{H}\otimes
I_{\alpha }^{H})\Delta _{1,\alpha }^{H}\mu _{\alpha }(g_{\alpha }\otimes
I_{\alpha })(c\otimes b) \\
&=&(I_{\alpha }\otimes \mu _{\alpha })(I_{\alpha }\otimes g_{\alpha
}^{-1}\otimes I_{\alpha })((\sigma _{\alpha }\otimes \sigma _{\alpha
^{-1}})\Delta _{\alpha ,\alpha ^{-1}}^{H}\otimes I_{\alpha }^{H}) \\
&&\Delta _{1,\alpha }^{H}\mu _{\alpha }(g_{\alpha }\otimes I_{\alpha
})(c\otimes b) \\
&=&(I_{\alpha }\otimes \mu _{\alpha })(I_{\alpha }\otimes g_{\alpha
}^{-1}\otimes I_{\alpha })\mu _{H_{\alpha }\otimes H_{\alpha ^{-1}}\otimes
H_{\alpha }} \\
&&(((\sigma _{\alpha }\otimes \sigma _{\alpha ^{-1}})\Delta _{\alpha ,\alpha
^{-1}}^{H}\otimes I_{\alpha }^{H})\Delta _{1,\alpha }^{H}g_{\alpha }\otimes
((\sigma _{\alpha }\otimes \sigma _{\alpha ^{-1}})\Delta _{\alpha ,\alpha
^{-1}}^{H}\otimes I_{\alpha }^{H})\Delta _{1,\alpha }^{H})(c\otimes b) \\
&=&(I_{\alpha }\otimes \mu _{\alpha })(I_{\alpha }\otimes g_{\alpha
}^{-1}\otimes I_{\alpha })\mu _{H_{\alpha }\otimes H_{\alpha ^{-1}}} \\
&&((I_{\alpha }\otimes I_{\alpha ^{-1}}\otimes g_{\alpha })(I_{\alpha
}\otimes \Delta _{\alpha ^{-1},\alpha }^{C})\Delta _{\alpha
,1}^{C}(c)\otimes 1\otimes 1\otimes b) \;\;\; \textrm{by Lemmas 3.8, 3.9}\\
&=&c_{1}^{\alpha }\otimes g_{\alpha }^{-1}(c_{21}^{\alpha ^{-1}})g_{\alpha
}(c_{22}^{\alpha })b=c_{1}^{\alpha }\otimes \epsilon (c_{2}^{1})b=\epsilon
(c_{2}^{1})c_{1}^{\alpha }\otimes b=c\otimes b.
\end{eqnarray*}%
\end{proof}

\begin{rem}
In any  Hopf  $\pi $ -coalgebra $H$, every $H_{\alpha }$ is
left $H_{1}$ -comodule by $\Delta _{1,\alpha }$. i.e., the following digrams
are commute
\begin{equation*}
\begin{tabular}{lllllllll}
& $H_{\alpha }$ & $%
\begin{tabular}{l}
$\Delta _{1,\alpha }$ \\
$\rightarrow $%
\end{tabular}%
$ & $H_{1}\otimes H_{\alpha }$ &  & $H_{\alpha }$ & $%
\begin{tabular}{l}
$\Delta _{1,\alpha }$ \\
$\rightarrow $%
\end{tabular}%
$ & $H_{1}\otimes H_{\alpha }$ &  \\
&  &  &  &  &  &  &  &  \\
$\Delta _{1,\alpha }$ & $\downarrow $ &  & $\downarrow \Delta _{1,1}\otimes
I_{\alpha }$ &  &  & $\sim \searrow $ & $\downarrow \epsilon ^{H}\otimes
I_{\alpha }$ &  \\
&  &  &  &  &  &  & $K\otimes H_{\alpha }$ &  \\
& $H_{1}\otimes H_{\alpha }$ &
\begin{tabular}{l}
$\longrightarrow $ \\
$I_{1}\otimes \Delta _{1,\alpha }$%
\end{tabular}
& $H_{1}\otimes H_{1}\otimes H_{\alpha }$ &  &  &  &  &
\end{tabular}%
\end{equation*}
\end{rem}

\begin{thm}
If $H$ have a left cosection, then $Ind(\rho )$ is
isomorphic to $V\otimes B=$ $\{(V\otimes B)_{\alpha }=V_{1}\otimes B_{\alpha
}\}_{\alpha \in \pi }$ as right $B$-module.
\end{thm}

\begin{proof}
Firstly, we`ll prove that $L_{1,\alpha }\eta _{\alpha
}g_{1}=(I_{1}\otimes \eta _{\alpha }g_{1})\Delta _{1,1}^{C}.$
\begin{eqnarray*}
L_{1,\alpha }\eta _{\alpha }g_{1}
&=&(\sigma _{1}\otimes I_{\alpha }^{H})\Delta _{1,\alpha }^{H}\eta _{\alpha
}g_{1} \\
&=&(\sigma _{1}\otimes I_{\alpha }^{H})(I_{1}\otimes \eta _{\alpha })\Delta
_{1,1}^{H}g_{1} \\
&=&(I_{1}\otimes \eta _{\alpha })(\sigma _{1}\otimes I_{1}^{H})\Delta
_{1,1}^{H}g_{1} \\
&=&(I_{1}\otimes \eta _{\alpha })L_{1,1}g_{1} \\
&=&(I_{1}\otimes \eta _{\alpha })(I_{1}\otimes g_{1})\Delta _{1,1}^{C} \\
&=&(I_{1}\otimes \eta _{\alpha }g_{1})\Delta _{1,1}^{C}.
\end{eqnarray*}

Now, we define $T_{\alpha }=(I_{1}\otimes \eta _{\alpha }g_{1})\rho
_{1,1}:V_{1}\rightarrow Ind(\rho )_{\alpha }.$ We`ll prove that $T_{\alpha
}(v_{1})\in Ind(\rho )_{\alpha }.$
\begin{eqnarray*}
(I_{1}\otimes L_{1,\alpha })T_{\alpha }(v_{1})
&=&(I_{1}\otimes L_{1,\alpha })(I_{1}\otimes \eta _{\alpha }g_{1})\rho
_{1,1}(v_{1}) \\
&=&(I_{1}\otimes L_{1,\alpha }\eta _{\alpha }g_{1})\rho _{1,1}(v_{1}) \\
&=&(I_{1}\otimes (I_{1}\otimes \eta _{\alpha }g_{1})\Delta _{1,1}^{C})\rho
_{1,1}(v_{1}) \\
&=&(I_{1}\otimes I_{1}\otimes \eta _{\alpha }g_{1})(I_{1}\otimes \Delta
_{1,1}^{C})\rho _{1,1}(v_{1}) \\
&=&(I_{1}\otimes I_{1}\otimes \eta _{\alpha }g_{1})(\rho _{1,1}\otimes
I_{1})\rho _{1,1}(v_{1}) \\
&=&(\rho _{1,1}\otimes I_{\alpha })(I_{1}\otimes \eta _{\alpha }g_{1})\rho
_{1,1}(v_{1}) \\
&=&(\rho _{1,1}\otimes I_{\alpha })T_{\alpha }(v_{1}).
\end{eqnarray*}

For $\alpha \in \pi $ we define $q_{\alpha }:V_{1}\otimes B_{\alpha
}\rightarrow Ind(\rho )_{\alpha }$ where
\begin{equation*}
q_{\alpha }(v\otimes b)=\lambda _{\alpha }(T_{\alpha }(v)\otimes
b)=v_{1}\otimes \eta _{\alpha }g_{1}(v_{2})b.
\end{equation*}%
Clear $q_{\alpha }(v\otimes b)\in Ind(\rho )_{\alpha }$ and $q_{\alpha }$ is
linear. We define
\begin{equation*}
q_{\alpha }^{-1}:Ind(\rho )_{\alpha }\rightarrow V_{1}\otimes B_{\alpha }
\end{equation*}%
where
\begin{eqnarray}
q_{\alpha }^{-1} &=&(I_{1}\otimes \mu _{\alpha }(\eta _{\alpha
}g_{1}^{-1}\otimes I_{\alpha }))(\rho _{1,1}\otimes I_{\alpha }) \label{3}\\
&=&(I_{1}\otimes \mu _{\alpha }(\eta _{\alpha }g_{1}^{-1}\otimes I_{\alpha
}))(I_{1}\otimes L_{1,\alpha }) \label{4}
\end{eqnarray}%

\begin{equation*}
i.e.,q_{\alpha }^{-1}(v\otimes h)=v_{1}\otimes \eta _{\alpha
}g_{1}^{-1}(v_{2})h=v\otimes \eta _{\alpha }g_{1}^{-1}\sigma
_{1}(h_{1}^{1})h_{2}^{\alpha }.
\end{equation*}
We`ll prove that $q_{\alpha }^{-1}(v\otimes h)\in V_{1}\otimes B_{\alpha }.$
\begin{eqnarray*}
(I_{1}\otimes L_{1,\alpha })q_{\alpha }^{-1}(v\otimes h)
&=&(I_{1}\otimes L_{1,\alpha })(I_{1}\otimes \mu _{\alpha }(\eta _{\alpha
}g_{1}^{-1}\otimes I_{\alpha }))(I_{1}\otimes L_{1,\alpha })(v\otimes h)\:\:\:\textrm{ by Equation \ref{4}} \\
&=&(I_{1}\otimes \mu _{C_{1}\otimes H_{\alpha }}(L_{1,\alpha }\eta _{\alpha
}g_{1}^{-1}\otimes L_{1,\alpha }))(I_{1}\otimes L_{1,\alpha })(v\otimes h) \\
&=&v\otimes \mu _{C_{1}\otimes H_{\alpha }}((I_{1}\otimes \eta _{\alpha
})L_{1,1}g_{1}^{-1}\otimes L_{1,\alpha })L_{1,\alpha }(h) \\
&=&v\otimes \mu _{C_{1}\otimes H_{\alpha }}((I_{1}\otimes \eta _{\alpha
})(S_{1}^{C}\otimes g_{1}^{-1})\tau \Delta _{1,1}^{H}\otimes L_{1,\alpha
})L_{1,\alpha }(h) \:\:\:\textrm{ by Lemma 3.4}\\
&=&v\otimes \mu _{C_{1}\otimes H_{\alpha }}((S_{1}^{C}\otimes \eta _{\alpha
}g_{1}^{-1})\tau \Delta _{1,1}^{H}\otimes (\sigma _{1}\otimes I_{\alpha
}^{H})\Delta _{1,\alpha }^{H}) (\sigma _{1}\otimes I_{\alpha }^{H})\Delta _{1,\alpha }^{H}(h) \\
&=&v\otimes \mu _{C_{1}\otimes H_{\alpha }}((S_{1}^{C}\otimes \eta _{\alpha
}g_{1}^{-1})\tau \otimes \sigma _{1}\otimes I_{\alpha }^{H})(\Delta _{1,1}^{H}\otimes \Delta _{1,\alpha }^{H})(\sigma _{1}\otimes
I_{\alpha }^{H})\Delta _{1,\alpha }^{H}(h) \\
&=&v\otimes \mu _{C_{1}\otimes H_{\alpha }}((S_{1}^{C}\otimes \eta _{\alpha
}g_{1}^{-1})\tau \otimes \sigma _{1}\otimes I_{\alpha }^{H}) \\
&&(\sigma _{1}\otimes \sigma _{1}\otimes I_{1}^{H}\otimes I_{\alpha
}^{H})(\Delta _{1,1}^{H}\otimes \Delta _{1,\alpha }^{H})\Delta _{1,\alpha
}^{H}(h) \\
&=&v\otimes S_{1}^{C}\sigma _{1}(h_{121}^{1})\sigma _{1}(h_{122}^{1})\otimes
\eta _{\alpha }g_{1}^{-1}\sigma _{1}(h_{11}^{1}))h_{2}^{\alpha } \\
&=&v\otimes \epsilon ^{C}\sigma _{1}(h_{12}^{1})1\otimes \eta _{\alpha
}g_{\alpha }^{-1}\sigma _{1}(h_{11}^{1}))h_{2}^{\alpha } \\
&=&v\otimes \epsilon ^{H}(h_{12}^{1})1\otimes \eta _{\alpha
}g_{1}^{-1}\sigma _{1}(h_{11}^{1}))h_{2}^{\alpha } \\
&=&v\otimes 1\otimes \eta _{\alpha }g_{1}^{-1}\sigma _{1}(\epsilon
^{H}(h_{12}^{1})h_{11}^{1}))h_{2}^{\alpha } \\
&=&v\otimes 1\otimes \eta _{\alpha }g_{1}^{-1}\sigma
_{1}(h_{1}^{1})h_{2}^{\alpha } \\
&=&v\otimes 1\otimes \mu _{\alpha }(\eta _{\alpha }g_{1}^{-1}\otimes
I_{\alpha }))L_{1,\alpha }(h)
\end{eqnarray*}%
Now, we will prove $q_{\alpha }q_{\alpha }^{-1}=I$ and $q_{\alpha }^{-1}q_{\alpha }=I$.
\begin{eqnarray*}
q_{\alpha }q_{\alpha }^{-1}(v\otimes h)
&=&(I_{1}\otimes \mu _{\alpha })(I_{1}\otimes \eta _{\alpha }g_{1}\otimes
I_{\alpha })(\rho _{1,1}\otimes I_{\alpha }) \\
&&(I_{1}\otimes \mu _{\alpha }(\eta _{\alpha }g_{1}^{-1}\otimes I_{\alpha
}))(\rho _{1,1}\otimes I_{\alpha })(v\otimes h) \:\:\:\textrm{ by Equation \ref{4}}\\
&=&(I_{1}\otimes \mu _{\alpha })(I_{1}\otimes \eta _{\alpha }g_{1}\otimes
I_{\alpha })(I_{1}\otimes I_{1}\otimes \mu _{\alpha }(\eta _{\alpha
}g_{1}^{-1}\otimes I_{\alpha })) \\
&&(\rho _{1,1}\otimes I_{1}\otimes I_{\alpha })(\rho _{1,1}\otimes I_{\alpha
})(v\otimes h) \\
&=&(I_{1}\otimes \mu _{\alpha })(I_{1}\otimes \eta _{\alpha }g_{1}\otimes
I_{\alpha })(I_{1}\otimes I_{1}\otimes \mu _{\alpha }(\eta _{\alpha
}g_{1}^{-1}\otimes I_{\alpha })) \\
&&((\rho _{1,1}\otimes I_{1})\rho _{1,1}\otimes I_{\alpha })(v\otimes h) \\
&=&(I_{1}\otimes \mu _{\alpha })(I_{1}\otimes \eta _{\alpha }g_{1}\otimes
I_{\alpha })(I_{1}\otimes I_{1}\otimes \mu _{\alpha }(\eta _{\alpha
}g_{1}^{-1}\otimes I_{\alpha })) \\
&&((I_{1}\otimes \Delta _{1,1}^{C}\otimes I_{\alpha })(\rho _{1,1}\otimes
I_{\alpha })(v\otimes h) \\
&=&v_{1}^{1}\otimes \eta _{\alpha }g_{1}(v_{21}^{1})\eta _{\alpha
}g_{1}^{-1}(v_{22}^{1})h\\
&=&v_{1}^{1}\otimes \eta _{\alpha }(g_{1}(v_{21}^{1})g_{1}^{-1}(v_{22}^{1}))h
\\
&=&v_{1}^{1}\otimes \eta _{\alpha }(\epsilon (v_{2}^{1})1)h \\
&=&\epsilon (v_{2}^{1})v_{1}^{1}\otimes \eta _{\alpha }(1)h \\
&=&v\otimes h
\end{eqnarray*}%
Also,
\begin{eqnarray*}
q_{\alpha }^{-1}q_{\alpha }(v\otimes b)
&=&(I_{1}\otimes \mu _{\alpha }(\eta _{\alpha }g_{1}^{-1}\otimes I_{\alpha
}))(I_{1}\otimes L_{1,\alpha })(I_{1}\otimes \mu _{\alpha }) \\
&&(I_{1}\otimes \eta _{\alpha }g_{1}\otimes I_{\alpha })(\rho _{1,1}\otimes
I_{\alpha })(v\otimes b) \:\:\:\textrm{ by Equation \ref{4}}\\
&=&(I_{1}\otimes \mu _{\alpha }(\eta _{\alpha }g_{1}^{-1}\otimes I_{\alpha
}))(I_{1}\otimes L_{1,\alpha })(v_{1}\otimes \eta _{\alpha
}g_{1}(v_{2}^{1})b) \\
&=&(I_{1}\otimes \mu _{\alpha }(\eta _{\alpha }g_{1}^{-1}\otimes I_{\alpha
}))(v_{1}\otimes L_{1,\alpha }(\eta _{\alpha }g_{1}(v_{2}^{1}))L_{1,\alpha
}(b)) \\
&=&(I_{1}\otimes \mu _{\alpha }(\eta _{\alpha }g_{1}^{-1}\otimes I_{\alpha
}))(v_{1}\otimes (I_{1}\otimes \eta _{\alpha })L_{1,1}g_{1}(v_{2}^{1})\cdot
(1\otimes b)) \\
&=&(I_{1}\otimes \mu _{\alpha }(\eta _{\alpha }g_{1}^{-1}\otimes I_{\alpha
}))(v_{1}\otimes (I_{1}\otimes \eta _{\alpha })(I_{1}\otimes g_{1})\Delta
_{1,1}^{C}(v_{2}^{1})\cdot (1\otimes b)) \\
&=&(I_{1}\otimes \mu _{\alpha }(\eta _{\alpha }g_{1}^{-1}\otimes I_{\alpha
}))(v_{1}\otimes v_{21}^{1}\otimes \eta _{\alpha }g_{1}(v_{22}^{1})b) \\
&=&v_{1}\otimes \eta _{\alpha }g_{1}^{-1}(v_{21}^{1})\eta _{\alpha
}g_{1}(v_{22}^{1})b \\
&=&v_{1}\otimes \eta _{\alpha }(g_{1}^{-1}(v_{21}^{1})g_{1}(v_{22}^{1}))b \\
&=&v_{1}\otimes \eta _{\alpha }(\epsilon ^{C}(v_{2}^{1})1)b \\
&=&\epsilon ^{C}(v_{2}^{1})v_{1}\otimes \eta _{\alpha }(1)b=v\otimes b.
\end{eqnarray*}

Now we`ll prove that $q_{\alpha }$ is module map for all $\alpha \in \pi ,$%
\begin{eqnarray*}
q_{\alpha }(I_{1}\otimes \mu _{\alpha })(v\otimes b\otimes k)
&=&q_{\alpha }(v\otimes bk) \\
&=&\lambda _{\alpha }(T_{\alpha }(v)\otimes bk) \\
&=&\lambda _{\alpha }((I_{1}\otimes \eta _{\alpha }g_{1})\rho
_{1,1}(v)\otimes bk)
\end{eqnarray*}%
\begin{eqnarray*}
&=&\lambda _{\alpha }(v_{1}\otimes \eta _{\alpha }g_{1}(v_{2})\otimes bk) \\
&=&v_{1}\otimes \eta _{\alpha }g_{1}(v_{2})bk \\
&=&v_{1}\otimes (\eta _{\alpha }g_{1}(v_{2})b)k \\
&=&\lambda _{\alpha }(v_{1}\otimes (\eta _{\alpha }g_{1}(v_{2})b)\otimes k)
\\
&=&\lambda _{\alpha }(q_{\alpha }\otimes I_{\alpha })(v\otimes b\otimes k).
\end{eqnarray*}

\end{proof}
Throughout this last part $(C,\sigma )$ is left $\pi -$ coisotropic quantum
subgroup of $H$ and $V=\{V_{\alpha }\}_{\alpha \in \pi }$ be a right $\pi -$%
comodule over $C$ and $G=\{G_{\alpha }\}_{\alpha \in \pi },$where
\begin{equation*}
G_{\alpha }=\{h\in H_{\alpha }:L_{1,\alpha }(h)=(\sigma _{1}\otimes
I_{\alpha }^{H})\Delta _{1,\alpha }^{H}(h)=\sigma _{1}(1)\otimes h\}.
\end{equation*}

\begin{lem}
$Ind(\rho )$ is right $G$-module
\end{lem}

\begin{proof}
Similar to Lemma 3.1.
\end{proof}

\begin{lem}
For $\alpha \in \pi ,u\in H_{1},c\in C_{\alpha ^{-1}}$
and $v\in \sigma _{\alpha ^{-1}}^{-1}(c)$ we have
\begin{eqnarray*}
&&(\sigma _{1}\otimes I_{\alpha }^{H})(\mu _{1}\otimes I_{\alpha
})(I_{1}\otimes \tau )(\Delta _{1,\alpha }^{H}g_{\alpha }^{-1}\otimes
I_{1})(c\otimes u) \\
&=&(I_{1}\otimes g_{\alpha }^{-1})(\sigma _{1}\otimes \sigma _{\alpha
^{-1}})(\mu _{1}(S_{1}^{H}\otimes I_{1})\otimes I_{\alpha
^{-1}})(I_{1}\otimes \tau )(\tau \Delta _{\alpha ^{-1},1}^{H}\otimes
I_{1})(v\otimes u)
\end{eqnarray*}
\end{lem}

\begin{proof}
 It can be proved with usual Hopf $\pi -$coalgebra
techniques.
\end{proof}

\begin{thm}
$H$ is isomorphic to $C\otimes G$ as vector space.
\end{thm}

\begin{proof}\mbox{}

\begin{enumerate}
\item For $h\in H_{\alpha },$ as in Lemma 3.6, $h=g_{\alpha }\sigma
_{\alpha }(h_{1}^{\alpha })g_{\alpha }^{-1}\sigma _{\alpha
^{-1}}(h_{21}^{\alpha ^{-1}})h_{22}^{\alpha }.$

\item For $h\in H_{1},\alpha \in \pi ,$ we`ll prove that
\begin{equation*}
\mu _{\alpha }(g_{\alpha }^{-1}\sigma _{\alpha ^{-1}}\otimes I_{\alpha
})\Delta _{\alpha ^{-1},\alpha }^{H}(h)\in G_{\alpha }.
\end{equation*}%
as follow%
\begin{eqnarray*}
L_{1,\alpha }\mu _{\alpha }(g_{\alpha }^{-1}\sigma _{\alpha ^{-1}}\otimes
I_{\alpha })\Delta _{\alpha ^{-1},\alpha }^{H}(h)
&=&(\sigma _{1}\otimes I_{\alpha }^{H})\Delta _{1,\alpha }^{H}\mu _{\alpha
}(g_{\alpha }^{-1}\sigma _{\alpha ^{-1}}\otimes I_{\alpha })\Delta _{\alpha
^{-1},\alpha }^{H}(h) \\
&=&(\sigma _{1}\otimes I_{\alpha }^{H})\mu _{H_{1}\otimes H_{\alpha
}}(\Delta _{1,\alpha }^{H}\otimes \Delta _{1,\alpha }^{H})(g_{\alpha
}^{-1}\sigma _{\alpha ^{-1}}\otimes I_{\alpha }^{H})\Delta _{\alpha
^{-1},\alpha }^{H}(h) \\
&=&(\sigma _{1}\mu _{H_{1}}\otimes \mu _{H_{\alpha }})(I\otimes \tau \otimes
I)(\Delta _{1,\alpha }^{H}\otimes \Delta _{1,\alpha }^{H})\\
&&(g_{\alpha
}^{-1}\sigma _{\alpha ^{-1}}\otimes I_{\alpha }^{H})\Delta _{\alpha
^{-1},\alpha }^{H}(h) \\
&=&(I_{1}\otimes \mu _{H_{\alpha }})((\sigma _{1}\mu _{H_{1}}\otimes
I_{\alpha })(I\otimes \tau )\otimes I_{\alpha }^{H})\\
&&(\Delta _{1,\alpha
}^{H}g_{\alpha }^{-1}\sigma _{\alpha ^{-1}}\otimes I_{1}\otimes I_{\alpha
}^{H}) (I_{\alpha ^{-1}}\otimes \Delta _{1,\alpha }^{H})\Delta _{\alpha
^{-1},\alpha }^{H}(h) \\
&=&(I_{1}\otimes \mu _{H_{\alpha }})((\sigma _{1}\mu _{H_{1}}\otimes
I_{\alpha })(I\otimes \tau )(\Delta _{1,\alpha }^{H}g_{\alpha }^{-1}\sigma
_{\alpha ^{-1}}\otimes I_{1})\otimes I_{\alpha }^{H}) \\
&&(I_{\alpha ^{-1}}\otimes \Delta _{1,\alpha }^{H})\Delta _{\alpha
^{-1},\alpha }^{H}(h) \\
&=&(I_{1}\otimes \mu _{H_{\alpha }})((I_{1}\otimes g_{\alpha }^{-1})(\sigma
_{1}\otimes \sigma _{\alpha ^{-1}})(\mu _{1}(S_{1}^{H}\otimes I_{1})\otimes
I_{\alpha ^{-1}}) \\
&&(I_{1}\otimes \tau )(\tau \Delta _{\alpha ^{-1},1}^{H}\otimes
I_{1})\otimes I_{\alpha }^{H})(I_{\alpha ^{-1}}\otimes \Delta _{1,\alpha
}^{H})\Delta _{\alpha ^{-1},\alpha }^{H}(h) \\
&=&(I_{1}\otimes \mu _{H_{\alpha }})((I_{1}\otimes g_{\alpha }^{-1})(\sigma
_{1}\otimes \sigma _{\alpha ^{-1}})(\mu _{1}(S_{1}^{H}\otimes I_{1})\otimes
I_{\alpha ^{-1}}) \\
&&(I_{1}\otimes \tau )(\tau \otimes I_{1})\otimes I_{\alpha }^{H})(\Delta
_{\alpha ^{-1},1}^{H}\otimes \Delta _{1,\alpha }^{H})\Delta _{\alpha
^{-1},\alpha }^{H}(h) \\
&=&(I_{1}\otimes \mu _{H_{\alpha }})[(\sigma _{1}\otimes g_{\alpha
}^{-1}\sigma _{\alpha ^{-1}})(\mu _{1}(S_{1}^{H}\otimes I_{1})\otimes
I_{\alpha ^{-1}}) \\
&&(I_{1}\otimes \tau )(\tau \otimes I_{1})\otimes I_{\alpha
}^{H}](h_{11}^{\alpha ^{-1}}\otimes h_{121}^{1}\otimes h_{122}^{1}\otimes
h_{2}^{\alpha }) \\
&=&\sigma _{1}(S_{1}^{H}(h_{121}^{1})h_{122}^{1})\otimes g_{\alpha
}^{-1}(\sigma _{\alpha ^{-1}}(h_{11}^{\alpha ^{-1}}))h_{2}^{\alpha } \\
&=&\sigma _{1}(\epsilon ^{H}(h_{12}^{1})1)\otimes g_{\alpha }^{-1}(\sigma
_{\alpha ^{-1}}(h_{11}^{\alpha ^{-1}}))h_{2}^{\alpha } \\
&=&\sigma _{1}(1)\otimes g_{\alpha }^{-1}(\sigma _{\alpha ^{-1}}(\epsilon
^{H}(h_{12}^{1})h_{11}^{\alpha ^{-1}}))h_{2}^{\alpha } \\
&=&\sigma _{1}(1)\otimes g_{\alpha }^{-1}(\sigma _{\alpha
^{-1}}(h_{1}^{\alpha ^{-1}}))h_{2}^{\alpha } \\
&=&\sigma _{1}(1)\otimes \mu _{\alpha }(g_{\alpha }^{-1}\sigma _{\alpha
^{-1}}\otimes I_{\alpha })\Delta _{\alpha ^{-1},\alpha }^{H}(h)
\end{eqnarray*}

\item We define $A_{\alpha }:C_{\alpha }\otimes G_{\alpha }\rightarrow
H_{\alpha }$ as follow%
\begin{equation*}
A_{\alpha }(c\otimes b)=\mu _{\alpha }(f_{\alpha }\otimes I_{\alpha
})(c\otimes b)=f_{\alpha }(c)b.
\end{equation*}
\end{enumerate}

Clear $A_{\alpha }$ is linear. We define $A_{\alpha }^{-1}:H_{\alpha
}\rightarrow C_{\alpha }\otimes K_{\alpha }$ as
\begin{equation*}
A_{\alpha }^{-1}(h)=(\sigma _{\alpha }\otimes \mu _{\alpha }(g_{\alpha
}^{-1}\sigma _{\alpha ^{-1}}\otimes I_{\alpha }))(\Delta _{\alpha ,\alpha
^{-1}}^{H}\otimes I_{\alpha })\Delta _{1,\alpha }^{H}(h).
\end{equation*}%
We`ll prove that for $\alpha \in \pi ,A_{\alpha }A_{\alpha
}^{-1}=I_{H_{\alpha }}$ and $A_{\alpha }^{-1}A_{\alpha }=I_{C_{\alpha
}\otimes B_{\alpha }}.$ Let $h\in H_{\alpha },c\otimes b\in C_{\alpha
}\otimes G_{\alpha },$%
\begin{eqnarray*}
A_{\alpha }A_{\alpha }^{-1}(h)
&=&A_{\alpha }((\sigma _{\alpha }\otimes \mu _{\alpha }(g_{\alpha
}^{-1}\sigma _{\alpha ^{-1}}\otimes I_{\alpha }))(\Delta _{\alpha ,\alpha
^{-1}}^{H}\otimes I_{\alpha })\Delta _{1,\alpha }^{H}(h)) \\
&=&A_{\alpha }((\sigma _{\alpha }\otimes \mu _{\alpha }(g_{\alpha
}^{-1}\sigma _{\alpha ^{-1}}\otimes I_{\alpha }))(I_{\alpha }\otimes \Delta
_{\alpha ^{-1},\alpha }^{H})\Delta _{\alpha ,1}^{H}(h)) \\
&=&A_{\alpha }(\sigma _{\alpha }(h_{1}^{\alpha })\otimes g_{\alpha
}^{-1}\sigma _{\alpha ^{-1}}(h_{21}^{\alpha ^{-1}})h_{22}^{\alpha }) \\
&=&g_{\alpha }\sigma _{\alpha }(h_{1}^{\alpha })g_{\alpha }^{-1}\sigma
_{\alpha ^{-1}}(h_{21}^{\alpha ^{-1}})h_{22}^{\alpha }=h.
\end{eqnarray*}%
Also,

\begin{eqnarray*}
A_{\alpha }^{-1}A_{\alpha }(c\otimes k)
&=&(\sigma _{\alpha }\otimes \mu _{\alpha }(g_{\alpha }^{-1}\sigma _{\alpha
^{-1}}\otimes I_{\alpha }))(\Delta _{\alpha ,\alpha ^{-1}}^{H}\otimes
I_{\alpha })\Delta _{1,\alpha }^{H}\mu _{\alpha }(g_{\alpha }\otimes
I_{\alpha })(c\otimes k) \\
&=&(I_{\alpha }\otimes \mu _{\alpha })(I_{\alpha }\otimes g_{\alpha
}^{-1}\otimes I_{\alpha })(\sigma _{\alpha }\otimes \sigma _{\alpha
^{-1}}\otimes I_{\alpha }^{H}) (\Delta _{\alpha ,\alpha ^{-1}}^{H}\otimes I_{\alpha }^{H})\Delta
_{1,\alpha }^{H}\mu _{\alpha }(g_{\alpha }\otimes I_{\alpha })(c\otimes k) \\
&=&(I_{\alpha }\otimes \mu _{\alpha })(I_{\alpha }\otimes g_{\alpha
}^{-1}\otimes I_{\alpha })((\sigma _{\alpha }\otimes \sigma _{\alpha
^{-1}})\Delta _{\alpha ,\alpha ^{-1}}^{H}\otimes I_{\alpha }^{H})\Delta
_{1,\alpha }^{H}\mu _{\alpha } (g_{\alpha }\otimes I_{\alpha })(c\otimes k) \\
&=&(I_{\alpha }\otimes \mu _{\alpha })(I_{\alpha }\otimes g_{\alpha
}^{-1}\otimes I_{\alpha })(\Delta _{\alpha ,\alpha ^{-1}}^{C}\sigma
_{1}\otimes I_{\alpha }^{H})\mu _{H_{\alpha }\otimes H_{\alpha ^{-1}}} (\Delta _{1,\alpha }^{H}g_{\alpha }\otimes \Delta _{1,\alpha
}^{H})(c\otimes k) \\
&=&(I_{\alpha }\otimes \mu _{\alpha })(I_{\alpha }\otimes g_{\alpha
}^{-1}\otimes I_{\alpha })(\Delta _{\alpha ,\alpha ^{-1}}^{C}\otimes
I_{\alpha }^{H})(\sigma _{1}\mu _{1}\otimes \mu _{\alpha })\\
&& (I\otimes \tau \otimes I)(\Delta _{1,\alpha }^{H}g_{\alpha }\otimes \Delta
_{1,\alpha }^{H})(c\otimes k) \\
&=&(I_{\alpha }\otimes \mu _{\alpha })(I_{\alpha }\otimes g_{\alpha
}^{-1}\otimes I_{\alpha })(\Delta _{\alpha ,\alpha ^{-1}}^{C}\otimes \mu
_{\alpha })\\
&&[(\sigma _{1}\mu _{1}\otimes I_{\alpha }^{H})(I\otimes \tau ) (\Delta _{1,\alpha }^{H}g_{\alpha }(c)\otimes k_{1}^{1}\otimes
k_{2}^{\alpha })] \\
&=&(I_{\alpha }\otimes \mu _{\alpha })(I_{\alpha }\otimes g_{\alpha
}^{-1}\otimes I_{\alpha })(\Delta _{\alpha ,\alpha ^{-1}}^{C}\otimes \mu
_{\alpha })\\
&& [(I_{1}\otimes g_{\alpha })(\sigma _{1}\otimes \sigma _{\alpha })(\mu
_{1}\otimes I_{\alpha })(I_{1}\otimes \tau )(\Delta _{1,\alpha }^{C}\otimes
I_{1})(v\otimes k_{1}^{1})\otimes k_{2}^{\alpha }] \\
&=&(I_{\alpha }\otimes \mu _{\alpha })(I_{\alpha }\otimes g_{\alpha
}^{-1}\otimes I_{\alpha })(\Delta _{\alpha ,\alpha ^{-1}}^{C}\otimes \mu
_{\alpha })[\sigma _{1}(v_{1}^{1}k_{1}^{1})\otimes g_{\alpha }\sigma
_{\alpha }(v_{2}^{\alpha })\otimes k_{2}^{\alpha }] \\
&=&(I_{\alpha }\otimes \mu _{\alpha })(I_{\alpha }\otimes g_{\alpha
}^{-1}\otimes I_{\alpha })(\Delta _{\alpha ,\alpha ^{-1}}^{C}\otimes \mu
_{\alpha }) [\omega _{1}(v_{1}^{1}\otimes \sigma _{1}(k_{1}^{1}))\otimes g_{\alpha
}\sigma _{\alpha }(v_{2}^{\alpha })\otimes k_{2}^{\alpha }] \\
&=&(I_{\alpha }\otimes \mu _{\alpha })(I_{\alpha }\otimes g_{\alpha
}^{-1}\otimes I_{\alpha })(\Delta _{\alpha ,\alpha ^{-1}}^{C}\otimes \mu
_{\alpha }) 
[\omega _{1}(v_{1}^{1}\otimes \sigma _{1}(1))\otimes g_{\alpha }\sigma
_{\alpha }(v_{2}^{\alpha })\otimes k] \\
&=&(I_{\alpha }\otimes \mu _{\alpha })(I_{\alpha }\otimes g_{\alpha
}^{-1}\otimes I_{\alpha })(\Delta _{\alpha ,\alpha ^{-1}}^{C}\otimes \mu
_{\alpha })[\sigma _{1}(v_{1}^{1})\otimes g_{\alpha }\sigma _{\alpha
}(v_{2}^{\alpha })\otimes k] \\
&=&(I_{\alpha }\otimes \mu _{\alpha })[(I_{\alpha }\otimes g_{\alpha
}^{-1})\Delta _{\alpha ,\alpha ^{-1}}^{C}\sigma _{1}(v_{1}^{1})\otimes
g_{\alpha }\sigma _{\alpha }(v_{2}^{\alpha })k]\\
&=&(I_{\alpha }\otimes \mu _{\alpha })[(I_{\alpha }\otimes g_{\alpha
}^{-1})(\sigma _{\alpha }\otimes \sigma _{\alpha ^{-1}})\Delta _{\alpha
,\alpha ^{-1}}^{H}(v_{1}^{1})\otimes g_{\alpha }\sigma _{\alpha
}(v_{2}^{\alpha })k] \\
&=&(I_{\alpha }\otimes \mu _{\alpha })[(I_{\alpha }\otimes g_{\alpha
}^{-1})(\sigma _{\alpha }(v_{1}^{\alpha })\otimes \sigma _{\alpha
^{-1}}(v_{2}^{\alpha ^{-1}})\otimes g_{\alpha }\sigma _{\alpha
}(v_{3}^{\alpha })k] \\
&=&\sigma _{\alpha }(v_{1}^{\alpha })\otimes g_{\alpha }^{-1}\sigma _{\alpha
^{-1}}(v_{2}^{\alpha ^{-1}})g_{\alpha }\sigma _{\alpha }(v_{3}^{\alpha })k \\
&=&\sigma _{\alpha }(v_{1}^{\alpha })\otimes \epsilon ^{C}\sigma
_{1}(v_{2}^{1})k \\
&=&\sigma _{\alpha }(\epsilon ^{H}(v_{2}^{1})v_{1}^{\alpha })\otimes k \\
&=&\sigma _{\alpha }(v)\otimes k \\
&=&c\otimes k.
\end{eqnarray*}

\end{proof}

\section{\large{ Coinduced representations of
Hopf group coalgebra}}

In this section, we study coinduced representation from left $\pi -$%
coisotropic quantum subgroup. We restrict our attention to finite
dimensional case of Hopf $\pi -$coalgebra, i.e., dim $H_{\alpha }=n^{\alpha
} $ $\prec \infty $ for all $\alpha \in \pi .$ Let us start now from left $%
\pi -$corepresentation $\rho =\{\rho _{\alpha ,\beta }\}_{\alpha ,\beta \pi }
$ of left $\pi -$coisotropic quantum subgroup $(C,\sigma )$ on $%
V=\{V_{\alpha }\}_{\alpha \in \pi }$. We define $W=\{W_{\alpha }\}_{\alpha
\in \pi },$ where $W_{\alpha }=\{F_{\alpha }\in Hom(V_{1},H_{\alpha }):$ $\
L_{1,\alpha }F_{\alpha }=(I_{C_{1}}\otimes F_{\alpha })\rho _{1,1}\}$ and $%
L_{\alpha ,\beta }=(\sigma _{\alpha }\otimes I_{\beta })\Delta _{\alpha
,\beta }$

\begin{lem}
For $\alpha ,\beta ,\gamma \in \pi ,$
\begin{equation*}
(L_{\alpha ,\beta }\otimes I_{\gamma })\Delta _{\alpha \beta ,\gamma
}=(I_{\alpha }\otimes \Delta _{\beta ,\gamma })L_{\alpha ,\beta \gamma }%
\end{equation*}
\end{lem}
\begin{proof}
Similar to Lemma 2.4 . 
\end{proof}
\begin{lem}
Suppose that $\{e_{i}^{\alpha }\}_{i=1}^{n^{\alpha }}$
and $\{g_{i}^{\alpha }\}_{i=1}^{n^{\alpha }}$ is its dual basis of $%
H_{\alpha }$ and $H_{\alpha }^{\ast }$ for all $\alpha \in \pi $. For $%
\alpha ,\beta ,\gamma \in \pi $, fix $F_{\alpha \beta \gamma }\in W_{\alpha
\beta \gamma }$, if we define the two maps $\xi _{1}:V_{1}\rightarrow
H_{\alpha }\otimes H_{\beta }\otimes H_{\gamma }$ and $\xi
_{2}:V_{1}\rightarrow H_{\alpha }\otimes H_{\beta }\otimes H_{\gamma }$ such
that
\begin{equation*}
\xi _{1}(v_{1})=\sum_{i=1}^{n^{\beta \gamma }}\sim (I_{\alpha }\otimes
g_{i}^{\beta \gamma })\Delta _{\alpha ,\beta \gamma }F_{\alpha \beta \gamma
}(v_{1})\otimes \Delta _{\beta ,\gamma }(e_{i}^{\beta \gamma })
\end{equation*}%
\begin{equation*}
\xi _{2}(v_{1})=\sum_{h=1}^{n^{\beta }}\sum_{l=1}^{n^{\gamma }}[\sim (\sim
\otimes I_{K})(I_{\alpha }\otimes (g_{h}^{\beta }\otimes g_{l}^{\gamma
})\Delta _{\beta ,\gamma }]\Delta _{\alpha ,\beta \gamma }F_{\alpha \beta
\gamma }(v_{1})\otimes e_{h}^{\beta }\otimes e_{l}^{\gamma }
\end{equation*}%
then $\xi _{1}=\xi _{2}$.
\end{lem}

\begin{proof}
We put

\begin{equation*}
F_{\alpha \beta \gamma }(v_{1})=h_{\alpha \beta \gamma
}=\sum_{j=1}^{n^{\alpha \beta \gamma }}\lambda _{j}e_{j}^{\alpha \beta
\gamma }
\end{equation*}%
\begin{equation*}
\Delta _{\beta ,\gamma }(e_{i}^{\beta \gamma })=\sum_{r=1}^{n^{\beta
}}\sum_{s=1}^{n^{\gamma }}\eta _{rs}^{i}e_{r}^{\beta }\otimes e_{s}^{\gamma }
\end{equation*}%
\begin{equation*}
\Delta _{\alpha ,\beta \gamma }(e_{i}^{\alpha \beta \gamma
})=\sum_{r=1}^{n^{\alpha }}\sum_{s=1}^{n^{\beta \gamma }}\theta
_{rs}^{i}e_{r}^{\alpha }\otimes e_{s}^{\beta \gamma }
\end{equation*}%
then we have

\begin{eqnarray*}
\Delta _{\alpha ,\beta \gamma }F_{\alpha \beta \gamma }(v_{1})
&=&\Delta _{\alpha ,\beta \gamma }(h_{\alpha \beta \gamma }) \\
&=&\sum_{j=1}^{n^{\gamma }}\lambda _{j}\Delta _{\alpha ,\beta \gamma
}(e_{j}^{\alpha \beta \gamma }) \\
&=&\sum_{j=1}^{n^{\gamma }}\sum_{r=1}^{n^{\alpha }}\sum_{s=1}^{n^{\beta
\gamma }}\lambda _{j}\theta _{rs}^{j}e_{r}^{\alpha }\otimes e_{s}^{\beta
\gamma }
\end{eqnarray*}%
imply that%
\begin{eqnarray*}
\sum_{i=1}^{n^{\beta \gamma }} &\sim &(I_{\alpha }\otimes g_{i}^{\beta
\gamma })\Delta _{\alpha ,\beta \gamma }F_{\alpha \beta \gamma
}(v_{1})\otimes \Delta _{\beta ,\gamma }(e_{i}^{\beta \gamma }) \\
&=&\sum_{i=1}^{n^{\beta \gamma }}[\sim (I_{\alpha }\otimes g_{i}^{\beta
\gamma })\sum_{j=1}^{n^{\gamma }}\sum_{r=1}^{n^{\alpha
}}\sum_{s=1}^{n^{\beta \gamma }}\lambda _{j}\theta _{rs}^{j}e_{r}^{\alpha
}\otimes e_{s}^{\beta \gamma }]\otimes \sum_{h=1}^{n^{\beta
}}\sum_{l=1}^{n^{\gamma }}\eta _{hl}^{i}e_{h}^{\beta }\otimes e_{l}^{\gamma }
\\
&=&\sum_{i=1}^{n^{\beta \gamma }}\sum_{j=1}^{n^{\gamma
}}\sum_{r=1}^{n^{\alpha }}\sum_{s=1}^{n^{\beta \gamma }}\lambda _{j}\theta
_{rs}^{j}g_{i}^{\beta \gamma }(e_{s}^{\beta \gamma })e_{r}^{\alpha }\otimes
\sum_{h=1}^{n^{\beta }}\sum_{l=1}^{n^{\gamma }}\eta _{hl}^{i}e_{h}^{\beta
}\otimes e_{l}^{\gamma } \\
&=&\sum_{i=1}^{n^{\beta \gamma }}\sum_{j=1}^{n^{\gamma
}}\sum_{r=1}^{n^{\alpha }}\lambda _{j}\theta _{ri}^{j}e_{r}^{\alpha }\otimes
\sum_{h=1}^{n^{\beta }}\sum_{l=1}^{n^{\gamma }}\eta _{hl}^{i}e_{h}^{\beta
}\otimes e_{l}^{\gamma } \\
&=&\sum_{i=1}^{n^{\beta \gamma }}\sum_{j=1}^{n^{\gamma
}}\sum_{r=1}^{n^{\alpha }}\sum_{h=1}^{n^{\beta }}\sum_{l=1}^{n^{\gamma
}}\eta _{hl}^{i}\lambda _{j}\theta _{ri}^{j}e_{r}^{\alpha }\otimes
e_{h}^{\beta }\otimes e_{l}^{\gamma }
\end{eqnarray*}%
Also
\begin{eqnarray*}
\sum_{h=1}^{n^{\beta }}\sum_{l=1}^{n^{\gamma }}[ &\sim &(\sim \otimes
I_{K})(I_{\alpha }\otimes (g_{h}^{\beta }\otimes g_{l}^{\gamma })\Delta
_{\beta ,\gamma }]\Delta _{\alpha ,\beta \gamma }F_{\alpha \beta \gamma
}(v_{1})\otimes e_{h}^{\beta }\otimes e_{l}^{\gamma } \\
&=&\sum_{h=1}^{n^{\beta }}\sum_{l=1}^{n^{\gamma }}[\sim (\sim \otimes
I_{K})(I_{\alpha }\otimes (g_{h}^{\beta }\otimes g_{l}^{\gamma })\Delta
_{\beta ,\gamma }] \\
&&\lbrack \sum_{j=1}^{n^{\gamma }}\sum_{r=1}^{n^{\alpha
}}\sum_{i=1}^{n^{\beta \gamma }}\lambda _{j}\theta _{ri}^{j}e_{r}^{\alpha
}\otimes e_{i}^{\beta \gamma }]\otimes e_{h}^{\beta }\otimes e_{l}^{\gamma }
\\
&=&\sum_{h=1}^{n^{\beta }}\sum_{l=1}^{n^{\gamma }}\sum_{j=1}^{n^{\gamma
}}\sum_{r=1}^{n^{\alpha }}\sum_{i=1}^{n^{\beta \gamma }}\lambda _{j}\theta
_{ri}^{j}\sim (\sim \otimes I_{K})[e_{r}^{\alpha }\otimes (g_{h}^{\beta
}\otimes g_{l}^{\gamma })\Delta _{\beta ,\gamma }(e_{i}^{\beta \gamma
})]\otimes e_{h}^{\beta }\otimes e_{l}^{\gamma } \\
&=&\sum_{h=1}^{n^{\beta }}\sum_{l=1}^{n^{\gamma }}\sum_{j=1}^{n^{\gamma
}}\sum_{r=1}^{n^{\alpha }}\sum_{i=1}^{n^{\beta \gamma }}\sum_{q=1}^{n^{\beta
}}\sum_{s=1}^{n^{\gamma }}\lambda _{j}\theta _{ri}^{j}\eta
_{qs}^{i}g_{h}^{\beta }(e_{q}^{\beta })g_{l}^{\gamma }(e_{s}^{\gamma
})e_{r}^{\alpha }\otimes e_{h}^{\beta }\otimes e_{l}^{\gamma } \\
&=&\sum_{h=1}^{n^{\beta }}\sum_{l=1}^{n^{\gamma }}\sum_{j=1}^{n^{\gamma
}}\sum_{r=1}^{n^{\alpha }}\sum_{i=1}^{n^{\beta \gamma }}\eta
_{hl}^{i}\lambda _{j}\theta _{ri}^{j}e_{r}^{\alpha }\otimes e_{h}^{\beta
}\otimes e_{l}^{\gamma }
\end{eqnarray*}%
therefore, we have
\begin{eqnarray*}
\sum_{i=1}^{n^{\beta \gamma }} &\sim &(I_{\alpha }\otimes g_{i}^{\beta
\gamma })\Delta _{\alpha ,\beta \gamma }F_{\alpha \beta \gamma
}(v_{1})\otimes \Delta _{\beta ,\gamma }(e_{i}^{\beta \gamma }) \\
&=&\sum_{h=1}^{n^{\beta }}\sum_{l=1}^{n^{\gamma }}[\sim (\sim \otimes
I_{K})(I_{\alpha }\otimes (g_{h}^{\beta }\otimes g_{l}^{\gamma })\Delta
_{\beta ,\gamma }]\Delta _{\alpha ,\beta \gamma }F_{\alpha \beta \gamma
}(v_{1})\otimes e_{h}^{\beta }\otimes e_{l}^{\gamma }
\end{eqnarray*}

then $\xi _{1}=\xi _{2}$
\end{proof}

\begin{thm}
$W=\{W_{\alpha }\}_{\alpha \in \pi }$ is right $\pi -$%
comodule over $H$ by $\Omega =\{\Omega _{\alpha ,\beta }\}_{\alpha ,\beta
\in \pi }$ where
\begin{equation*}
\Omega _{\alpha ,\beta }:W_{\alpha \beta }\rightarrow W_{\alpha }\otimes
H_{\beta }\text{ as }\Omega _{\alpha ,\beta }(F_{\alpha \beta })=\sim
(I_{\alpha }\otimes g^{\beta })\Delta _{\alpha ,\beta }F_{\alpha \beta
}\otimes e^{\beta }.\newline
\end{equation*}
\end{thm}

\begin{proof}
Firstly, we ${^{\prime }}$ll prove that $\sim (I_{\alpha
}\otimes g^{\beta })\Delta _{\alpha ,\beta }F_{\alpha \beta }\in W_{\alpha }$.
\begin{eqnarray*}
L_{1,\alpha }[\sim (I_{\alpha }\otimes g^{\beta })\Delta _{\alpha ,\beta
}F_{\alpha \beta }]
&=&\sim (L_{1,\alpha }\otimes g^{\beta })\Delta _{\alpha ,\beta }F_{\alpha
\beta } \\
&=&(I_{C_{1}}\otimes \sim (I_{\alpha }\otimes g^{\beta }))(L_{1,\alpha
}\otimes I_{\beta })\Delta _{\alpha ,\beta }F_{\alpha \beta } \\
&=&(I_{C_{1}}\otimes \sim (I_{\alpha }\otimes g^{\beta }))(I_{C_{1}}\otimes
\Delta _{\alpha ,\beta })L_{1,\alpha \beta }F_{\alpha \beta } \\
&=&(I_{C_{1}}\otimes \sim (I_{\alpha }\otimes g^{\beta }))(I_{C_{1}}\otimes
\Delta _{\alpha ,\beta })(I_{C_{1}}\otimes F_{\alpha \beta })\rho _{1,1} \\
&=&(I_{C_{1}}\otimes \sim (I_{\alpha }\otimes g^{\beta })\Delta _{\alpha
,\beta }F_{\alpha \beta })\rho _{1,1} 
\end{eqnarray*}%

Now, we will prove that $\Omega $ is coaction on $W,$ i.e., the
two diagrams \newline

\begin{tabular}{lll}
$\ \ \ \ \ W_{\alpha \beta \gamma }$ & $%
\begin{tabular}{l}
$\Omega _{\alpha \beta ,\gamma }$ \\
$\rightarrow $%
\end{tabular}%
$ & $W_{\alpha \beta }\otimes H_{\gamma }$ \\
&  &  \\
$\Omega _{\alpha ,\beta \gamma }\downarrow $ &  & $\downarrow (\Omega
_{\alpha ,\beta }\otimes I_{\gamma })$ \\
&  &  \\
$W_{\alpha }\otimes H_{\beta \gamma }$ &
\begin{tabular}{l}
$\longrightarrow $ \\
$I_{W_{\alpha }}\otimes \Delta _{\beta ,\gamma }$%
\end{tabular}
& $W_{\alpha }\otimes H_{\beta }\otimes H_{\gamma }$%
\end{tabular}%
, \ \ \
\begin{tabular}{lll}
$W_{\alpha }$ & $%
\begin{tabular}{l}
$\Omega _{\alpha ,1}$ \\
$\rightarrow $%
\end{tabular}%
$ & $W_{\alpha }\otimes H_{1}$ \\
&  &  \\
& $\sim $ $\searrow $ & $\downarrow (I_{W_{\alpha }}\otimes \epsilon )$ \\
&  &  \\
&  & $W_{\alpha }\otimes K$%
\end{tabular}

are commutes.

\begin{enumerate}
\item
\begin{eqnarray*}
(\Omega _{\alpha ,\beta }\otimes I_{\gamma })\Omega _{\alpha \beta ,\gamma
}F_{\alpha \beta \gamma }
&=&(\Omega _{\alpha ,\beta }\otimes I_{\gamma })(\sim (I_{\alpha \beta
}\otimes g^{\gamma })\Delta _{\alpha \beta ,\gamma }F_{\alpha \beta \gamma
}\otimes e^{\gamma }) \\
&=&\Omega _{\alpha ,\beta }(\sim (I_{\alpha \beta }\otimes g^{\gamma
})\Delta _{\alpha \beta ,\gamma }F_{\alpha \beta \gamma })\otimes e^{\gamma }
\\
&=&\sim (I_{\alpha }\otimes g^{\beta })\Delta _{\alpha ,\beta }\sim
(I_{\alpha \beta }\otimes g^{\gamma })\Delta _{\alpha \beta ,\gamma
}F_{\alpha \beta \gamma }\otimes e^{\beta }\otimes e^{\gamma } \\
&=&\sim (\sim \otimes I_{K})((I_{\alpha }\otimes g^{\beta })\Delta _{\alpha
,\beta }\otimes g^{\gamma })\Delta _{\alpha \beta ,\gamma }F_{\alpha \beta
\gamma }\otimes e^{\beta }\otimes e^{\gamma } \\
&=&\sim (\sim \otimes I_{K})(I_{\alpha }\otimes g^{\beta }\otimes g^{\gamma
})(\Delta _{\alpha ,\beta }\otimes I_{\gamma })\Delta _{\alpha \beta ,\gamma
}F_{\alpha \beta \gamma }\otimes e^{\beta }\otimes e^{\gamma } \\
&=&\sim (\sim \otimes I_{K})(I_{\alpha }\otimes g^{\beta }\otimes g^{\gamma
})(I_{\alpha }\otimes \Delta _{\beta ,\gamma })\Delta _{\alpha ,\beta \gamma
}F_{\alpha \beta \gamma }\otimes e^{\beta }\otimes e^{\gamma } \\
&=&\sim (\sim \otimes I_{K})(I_{\alpha }\otimes (g^{\beta }\otimes g^{\gamma
})\Delta _{\beta ,\gamma })\Delta _{\alpha ,\beta \gamma }F_{\alpha \beta
\gamma }\otimes e^{\beta }\otimes e^{\gamma }= \xi_1
\end{eqnarray*}%
Also
\begin{eqnarray*}
(I_{W_{\alpha }}\otimes \Delta _{\beta ,\gamma })\Omega _{\alpha ,\beta
\gamma }F_{\alpha \beta \gamma }
&=&(I_{W_{\alpha }}\otimes \Delta _{\beta ,\gamma })(\sim (I_{\alpha
}\otimes g^{\beta \gamma })\Delta _{\alpha ,\beta \gamma }F_{\alpha \beta
\gamma }\otimes e^{\beta \gamma }) \\
&=&\sim (I_{W_{\alpha }}\otimes g^{\beta \gamma })\Delta _{\alpha ,\beta
\gamma }F_{\alpha \beta \gamma }\otimes \Delta _{\beta ,\gamma }(e^{\beta
\gamma })=\xi_2.
\end{eqnarray*}%
from Lemma 4.2. the first digram is commut

\item
\begin{eqnarray*}
(I_{W_{\alpha }}\otimes \epsilon )\Omega _{\alpha ,1}F_{\alpha }
&=&(I_{W_{\alpha }}\otimes \epsilon )[\sim (I_{\alpha }\otimes g^{1})\Delta
_{\alpha ,1}F_{\alpha }\otimes e^{1}] \\
&=&\sim (I_{W_{\alpha }}\otimes g^{1})\Delta _{\alpha ,1}F_{\alpha }\otimes
\epsilon (e^{1}) \\
&=&\sim (I_{W_{\alpha }}\otimes \epsilon (e^{1})g^{1})\Delta _{\alpha
,1}F_{\alpha }\otimes 1_{K} \\
&=&\sim (I_{W_{\alpha }}\otimes \epsilon )\Delta _{\alpha ,1}F_{\alpha
}\otimes 1_{K} \\
&=&F_{\alpha }\otimes 1_{K}.
\end{eqnarray*}
\end{enumerate}
\end{proof}

Now, we constract an induced and coinduced representation from subHopf $\pi
- $coalgebra.

\begin{thm}
Let $H$ be a finite dimentional Hopf $\pi $-coalgebra
and $A=\{A_{\alpha }\}_{\alpha \in \pi }$ be an isolated subHopf $\pi $%
-coalgebra of $H.$ If $V=\{V_{\alpha }\}_{\alpha \in \pi }$ is left $\pi $%
-comodule over $A$ by $\rho =\{\rho _{\alpha ,\beta }\}_{\alpha ,\beta \in
\pi },$ then we can construct an induced and coinduced representation over $%
H $.
\end{thm}

\begin{proof}
 Since $A=\{A_{\alpha }\}_{\alpha \in \pi }$ is an isolated
subHopf $\pi $-coalgebra of $H$, there exist a family $I=\{I_{\alpha
}\}_{\alpha \in \pi }$ of Hopf $\pi $-coideal of $H$ such that $H_{\alpha
}=A_{\alpha }\oplus I_{\alpha }$ for all $\alpha \in \pi .$ We$^{,}$ll prove
that $(A,\sigma )$ is left $\pi -$coisotropic quantum subgroup of $H$ where $%
\sigma =\{\sigma _{\alpha }\}_{\alpha \in \pi }$ and $\sigma _{\alpha
}:H_{\alpha }\rightarrow A_{\alpha }$ as $\sigma _{\alpha }(m_{\alpha
}+i_{\alpha })=m_{\alpha }.$\\

Clear, $A$ is $\pi $-coalgebra. We ${^{\prime }}$ll prove that $A_{\alpha }$
is left $H_{\alpha }$-module for all $\alpha \in \pi .$ We define $\Phi
_{\alpha }:H_{\alpha }\otimes A_{\alpha }\rightarrow A_{\alpha }$ as follow $%
\Phi _{\alpha }((m_{\alpha }+i_{\alpha })\otimes a_{\alpha })=m_{\alpha
}a_{\alpha }$. Then
\begin{eqnarray*}
\Phi _{\alpha }(\mu _{\alpha }\otimes I_{\alpha })((m_{\alpha }+i_{\alpha
})\otimes (n_{\alpha }+j_{\alpha })\otimes a_{\alpha })
&=&\Phi _{\alpha }((m_{\alpha }n_{\alpha }+m_{\alpha }j_{\alpha }+i_{\alpha
}n_{\alpha }+i_{\alpha }j_{\alpha })\otimes a_{\alpha })\\&=&m_{\alpha
}n_{\alpha }a_{\alpha }
\end{eqnarray*}%
and
\begin{eqnarray*}
\Phi _{\alpha }(I_{\alpha }\otimes \Phi _{\alpha })((m_{\alpha }+i_{\alpha
})\otimes (n_{\alpha }+j_{\alpha })\otimes a_{\alpha })
&=&\Phi _{\alpha }((m_{\alpha }+i_{\alpha })\otimes n_{\alpha }a_{\alpha
})\\
&=&m_{\alpha }n_{\alpha }a_{\alpha }
\end{eqnarray*}%
where $m_{\alpha }j_{\alpha }+i_{\alpha }n_{\alpha }+i_{\alpha }j_{\alpha
}\in I_{\alpha }$\\
Also, we have
\begin{eqnarray*}
\Phi _{\alpha }(\eta _{\alpha }\otimes I_{\alpha })(k\otimes a_{\alpha })
&=&\Phi _{\alpha }((\eta _{\alpha }(k)+0)\otimes a_{\alpha }) \\
&=&\eta _{\alpha }(k)a_{\alpha }=ka_{\alpha }.
\end{eqnarray*}

Now, we will prove that $\sigma _{\alpha }$ is $\pi -$coalgebra map
\begin{eqnarray*}
(\sigma _{\alpha }\otimes \sigma _{\beta })\Delta _{\alpha ,\beta
}(h_{\alpha \beta })
&=&(\sigma _{\alpha }\otimes \sigma _{\beta })\Delta _{\alpha ,\beta
}(m_{\alpha \beta }+i_{\alpha \beta }) \\
&=&(\sigma _{\alpha }\otimes \sigma _{\beta })(\Delta _{\alpha ,\beta
}(m_{\alpha \beta })+\Delta _{\alpha ,\beta }(i_{\alpha \beta })) \\
&=&\Delta _{\alpha ,\beta }(m_{\alpha \beta })
\end{eqnarray*}%
and
\begin{eqnarray*}
\Delta _{\alpha ,\beta }\sigma _{\alpha \beta }(h_{\alpha \beta })
&=&\Delta _{\alpha ,\beta }\sigma _{\alpha \beta }(m_{\alpha \beta
}+i_{\alpha \beta }) \\
&=&\Delta _{\alpha ,\beta }(m_{\alpha \beta }).
\end{eqnarray*}
We will prove that $\sigma _{\alpha }$ is left module map
\begin{eqnarray*}
\Phi _{\alpha }(I_{\alpha }\otimes \sigma _{\alpha })(h_{\alpha }\otimes
k_{\alpha })
&=&\Phi _{\alpha }(I_{\alpha }\otimes \sigma _{\alpha })((m_{\alpha
}+i_{\alpha })\otimes (n_{\alpha }+j_{\alpha })) \\
&=&\Phi _{\alpha }((m_{\alpha }+i_{\alpha })\otimes n_{\alpha })=m_{\alpha
}n_{\alpha }
\end{eqnarray*}%
and
\begin{eqnarray*}
\sigma _{\alpha }\mu _{\alpha }(h_{\alpha }\otimes k_{\alpha })
&=&\sigma _{\alpha }(m_{\alpha }n_{\alpha }+m_{\alpha }j_{\alpha }+i_{\alpha
}n_{\alpha }+i_{\alpha }j_{\alpha })=m_{\alpha }n_{\alpha }
\end{eqnarray*}

Therefore $(A,\sigma )$ is left $\pi -$coisotropic quantum subgroup of $H.$
Since $V$ is left $\pi -$comodule over $A$, then we can construct an induced
and coinduced representation over $H.$
\end{proof}

\end{document}